\documentclass[a4paper,11pt]{article}
\usepackage{color}
\usepackage{float}
\usepackage{RR}
\usepackage{hyperref}
\usepackage{pstricks,pst-node,pst-text,pst-3d}
\usepackage{rotate,float,pst-plot}
\usepackage{latexsym}
\usepackage{amssymb}
\usepackage[latin1]{inputenc}
\usepackage{fancyhdr}
\usepackage{amsmath}
\usepackage{amsthm}
\usepackage{a4wide}
\usepackage[amssymb]{SIunits}
\RRdate{Juillet 2008} 

\RRauthor{
Julien Diaz
  \thanks[sfn]{EPI Magique-3D, Centre de Recherche
  Inria Bordeaux Sud-Ouest}
\thanks[sf1]{Laboratoire de Mathématiques et de
  leurs Applications, CNRS UMR-5142, Université de Pau et des Pays de
  l'Adour --  Bâtiment IPRA, avenue de
  l'Université  -- BP 1155-64013 PAU CEDEX}%
  \and
Abdelaâziz Ezziani
\thanksref{sf1}
\thanksref{sfn}
}
\authorhead{Diaz \& Ezziani}
\RRtitle{Solution analytique pour la propagation
  d'ondes en milieu poroélastique stratifié. Partie~II~: en dimension 3}
\RRetitle{Analytical Solution for Wave Propagation in Stratified Poroelastic Medium. Part~II: the 3D Case} 
\titlehead{Analytical solution}
\RRresume{Nous nous intéressons à la modélisation
  de la propagation d'ondes dans les milieux
  infinis bicouches poroélastiques.  Nous
  considérons ici le modèle bi-phasique de Biot. 
Cette seconde partie
  est consacrée au calcul de la solution
  analytique en dimension trois à l'aide de la
  technique de Cagniard-De Hoop. }
\RRabstract{We are interested in the modeling of wave propagation in poroelastic media.
  We consider the biphasic Biot's model in an infinite bilayered medium,
  with a plane interface.
  We adopt the Cagniard-De Hoop's technique. This report is devoted to the
  calculation of analytical solution in three dimension.} 
\RRmotcle{Modèle de Biot, ondes poroélastiques, solution analytique,
technique de Cagniard de Hoop.
 } 
\RRmotcle{Modèle de Biot, ondes poroélastiques, solution analytique,
  technique de Cagniard de Hoop.}
\RRkeyword{Biot's model, poroelastic waves,
  analytical solution, Cagniard-De Hoop's
  technique.} 
\RRprojet{Magique-3D} 
\RRtheme{\THNum} 
 \URFuturs 
\RCBordeaux
\usepackage{ncacousporo}
\RRNo{6596}
\begin{document}
\makeRR
\section*{Introduction}
Many seismic materials  cannot only be considered  as solid materials. 
They are often porous media, i.e. media made of a solid fully
saturated with a fluid: there are solid media perforated by a
multitude of small holes (called pores) filled with a fluid. 
It is in particular often the case of the oil 
reservoirs. It is clear that the analysis of results by seismic 
methods of the exploration of such media must take to account the fact that a
wave being propagated in such a medium meets a succession of phases
solid and fluid: we speak about poroelastic media, and the more
commonly used model is the Biot's model \cite{biot1,biot2,biot3}. \\\\
 When the wavelength is large in comparison
with the size of the pores, rather than  regarding such a medium as
an heterogeneous medium, it is legitimate to use, at least locally,
the theory of homogenization~\cite{burridge,hom_por}. This leads
to the Biot's model~\cite{biot1,biot2,biot3} which involves  as
unknown not only the displacement field in the solid but
also the displacement field in the fluid. The principal characteristic of this
model is that in addition to  the classical P and S waves in a solid
one observes a P ``slow'' wave, which we could also call a  ``fluid''
wave: the denomination ``slow wave'' refers to the fact that in
practical applications, it is slower (and probably much slower) than
the other two waves.\\\\
The computation of analytical solutions for wave
propagation in poroelastic media is of high importance for the
validation of numerical computational codes or for
a better understanding of the
reflexion/transmission properties of the media.
Cagniard-de Hoop method~\cite{Cag,DH} is a
useful tool to obtain such solutions and permits
to compute each type of waves
(P wave, S wave, head wave...) independently.
Although
it was originally dedicated to the solution to
elastodynamic wave propagation, it can be applied
to any transient wave propagation problem in
stratified medium. However, as far as we know,
few works have been dedicated to the
application of this method to poroelastic medium, especially in three dimensions.

In order to validate computational codes of wave
propagation in poroelastic media, we have
implemented the codes Gar6more 2D~\cite{Gar6} and Gar6more 3D~\cite{Gar63d}
which provide the complete solution (reflected
and transmitted waves) of the propagation of wave
in stratified 2D or 3D media composed of
acoustic/acoustic, acoustic/elastic,
acoustic/poroelastic or poroelastic/poroelastic
layers. The codes are freely downloadable at\\\centerline{\url{http://www.spice-rtn.org/library/software/Garcimore2D}} and 
\\\centerline{\url{http://www.spice-rtn.org/library/software/Gar6more3D}.}
\\\\
We will focus in this paper on the 3D poroelastic
case, the two dimensional and the acoustic/poroelastic
cases are detailed in~\cite{RAP_DE6509,RAP2,RAP3}. The
outline of the paper is as follows: we first
present the model problem we want to solve and
derive the Green problem from it (section 1).
Then we present the analytical solution to the wave
propagation problem in a stratified 2D medium composed of
an acoustic and a poroelastic layer (section 2).
Finally we illustrate our results through numerical applications (section 3).
\section{The model problem}
\begin{figure}[htbp]
  \centering
\setlength{\unitlength}{.9cm}
  \begin{picture}(6,6.5)(2,0)
\psset{xunit=1.5cm}
\psframe[linestyle=dotted](0,.25)(6,5.75)
\psline(0,3)(6,3)
\put(5,4.3){$\Omega^+$}
\put(5,1.){$\Omega^-$}
\put(3.5,5.2){First Layer}
\put(3.5,2){Second Layer}
\put(-1.2,3.2){$z=0$}
  \end{picture}
\caption{Configuration of the study}
\label{fig:interf_plan}
\end{figure}
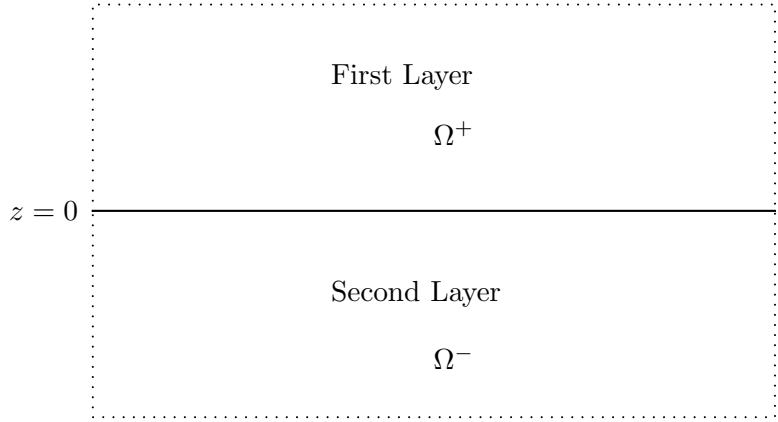

We consider an infinite two dimensional medium
($\Om=\R^3$) composed of two homogeneous poroelastic layers  
$\Omega^+=\R^2\times]-\infty,0]$ and $\Omega^-=\R^2\times[0,+\infty[$ separated by
an horizontal interface $\Gamma$ (see
Fig.~\ref{fig:interf_plan}). We first describe the
equations in the two layers
(\S\ref{sec:equation-acoustics})  and the transmission
conditions on the interface $\Gamma$
(\S\ref{sec:transm-cond}), then we present
the Green problem from which we  compute the
analytical solution (\S\ref{sec:greens-problem}).
\subsection{Poroelastic equations}
\label{sec:equation-acoustics}
We consider the second-order
formulation of the poroelastic equations~\cite{biot1,biot2,biot3}: 
\begin{equation}\label{eq:biot}
\left\{
\begin{array}{lll}
\dsp\ro\,\ddU+\ro_f\,\ddW-\Nabla\cdot\Si=\gr F_u,&&\mbox{in }
\Om\times]0,T],\\[8pt]
\dsp\ro_f\,\ddU+\ro_w\,\ddW+\frac{1}{\mat K}\,\gr \dW+\nabla P=\gr F_w,
&&\mbox{in }\Om\times]0,T],\\[12pt]
\dsp \Si=\la\nabla\cdot \gr U_s\,\gr I_3+2\mu\vare(\gr U_s)-\be\, P\,\gr I_3, &&\mbox{in }\Om\times]0,T],
\\[8pt]
\dsp\frac{1}{m}\,P+\be\,\nabla\cdot\gr U_s+\nabla\cdot \gr
W=F_p,&&\mbox{in }\Om\times]0,T],\\[12pt]
\gr U_s(x,0)=0,\,\gr W(x,0)=0, &&
\mbox{in }\Om,\\[12pt]
\dU(x,0)=0,\,\dW(x,0)=0, &&
\mbox{in }\Om,
\end{array}
\right.
\end{equation}
with
$$
(\Nabla\cdot\Si)_i=\sum_{j=1}^3\frac{\p\Si_{ij}}{\p
  x_j}\;\;\forall\,i=1,3.\;\mbox{ As usual }\gr{I}_3 \mbox{ is the identity matrix of }\mat M_2(\RR),
$$
and $\vare(\gr U_s)$ is the solid strain tensor defined by: 
$$
\vare_{ij}(\gr U)=\fr\left(\frac{\p U_i}{\p x_j}+\frac{\p U_j}{\p x_i}\right).$$
In (\ref{eq:biot}), the unknowns are:\\
\begin{itemize}
\item $\gr U_s$ the displacement field of solid particles;\\
\item $\gr W=\phi(\gr U_f-\gr U_s)$, the relative displacement, $\gr
  U_f$ being the displacement field of fluid particles and $\phi$ 
the porosity;\\
\item $P$, the fluid pressure;\\
\item $\Si$, the solid stress tensor.\\
\end{itemize}
The parameters describing the physical properties
of the medium are given by:\\
\begin{itemize}
\item $\ro=\phi\,\ro_f+(1-\phi)\ro_s$ is the overall density of the saturated medium, with $\ro_s$ the density of the solid and  $\ro_f$ the density of the fluid;\\
\item $\ro_w=a\ro_f/\phi$, where $a$ is the tortuosity of the solid matrix;\\
\item $\mat K=\kappa/\eta$, where $\kappa$ is the permeability of the solid matrix and $\eta\-$ is the viscosity of the fluid;\\
\item $m$ and $\beta$ are positive physical coefficients: 
$\be=1-K_b/K_s$ \\
and $m=\left[\phi/K_f+(\be-\phi)/K_s\right]^{-1}$, where 
$K_s$ is the bulk modulus of the solid, $K_f$ is the bulk modulus of the fluid 
and $K_b$ is the frame bulk modulus;\\
\item  $\mu$ is the frame shear modulus,
and $\la=K_b-2\mu/3$ is the Lam\'e constant.\\
\item $\gr F_u$, $\gr F_w$ and $F_p$ are the force densities.
\end{itemize}
To simplify this study, we consider only the case of a compression source 
$$\gr F_u(x,y,t)=f_u\nabla(\delta_x\,\delta_y\,\delta_{z-h})\,f(t) \hbox{ and } \gr F_w(x,y,t)=f_w\nabla(\delta_x\,\delta_y\,\delta_{z-h})\,f(t)$$ and a pressure source $F_p=f_p\delta_x\,\delta_y\,\delta_{z-h}\,f(t)$, where $f_u$, $f_w$ and $f_p$ are constant and $f$ is a regular source function in time. We can generalize this 
approach for other types of punctual sources such as for instance $$\gr F_u=f_u\nabla\times(\delta_x\,\delta_y\,\delta_{z-h}\,\gr v)\,f(t)\hbox{ and }\gr F_w=f_w\nabla\times(\delta_x\,\delta_y\,\delta_{z-h}\,\gr v)\,f(t)$$
where $\gr v$ is a vector in $\RR^3$.
\subsection{Transmission conditions}
\label{sec:transm-cond}
Let $\gr n$ be the unitary normal vector of
$\Gamma$ outwardly directed to $\Om^-$. The transmission 
conditions on the interface $\Ga$ between the two poroelastic medium  
are~\cite{carcione}:
\begin{equation}\label{eq:condi_trans}
\left\{
\begin{array}{l}
\gr U_s^+=\gr U_s^-,\\[4pt]
\gr W^+ \cdot \gr n=\gr W^-\cdot\gr n,\\[4pt]
P^+=P^-,\\[4pt]
\Si^+\,\gr n=\Si^-\,\gr n.
\end{array}
\right.
\end{equation}
\subsection{The Green problem}
\label{sec:greens-problem}
We won't compute directly the solution to
(\ref{eq:condi_trans})
but the solution to the following Green problem:
\begin{subequations}\label{biot_decomp}
\begin{eqnarray}
\dsp\ro^\pm\,\ddu^\pm+\ro_f^\pm\,\ddw^\pm-\Nabla\cdot\si^\pm=
f_u\,\nabla(\delta_x\,\delta_y\,\delta_{z-h})\delta_t,&&\mbox{in }
\Om^\pm\times]0,T],\label{eq:biot1}\\[8pt]
\dsp\ro_f^\pm\,\ddu^\pm+\ro_w^\pm\,\ddw^\pm+\frac{1}{\mat K^\pm}\,\gr \dw^\pm+\nabla p^\pm=f_w\,\nabla(\delta_x\,\delta_y\,\delta_{z-h})\delta_t,
&&\mbox{in }\Om^\pm\times]0,T],\label{eq:biot2}\\[8pt]
\dsp \si^\pm=\la^\pm\nabla\cdot \gr u^\pm_s\,\gr I_3+2\mu^\pm\vare(\gr u^\pm_s)- \be^\pm\, p^\pm\,\gr I_3, &&\mbox{in }\Om^\pm\times]0,T],
\label{acousporo:eq:6}\\[8pt]
\dsp\frac{1}{m^\pm}\,p^\pm+\be^\pm\,\nabla\cdot\gr u^\pm_s+\nabla\cdot \gr
w^\pm=f_p\,\delta_x\,\delta_y\,\delta_{z-h}\,f(t),&&\mbox{in }\Om^\pm\times]0,T],\label{acousporo:eq:8}\\[16pt]
\gr u_s^-=\gr u^+,&&\mbox{on }\Gamma\times]0,T]\\[4pt]
\gr w^-\cdot\gr n=\gr w^+\cdot\gr n,&&\mbox{on }\Gamma\times]0,T]\\[4pt]
p^-=p^+,&&\mbox{on }\Gamma\times]0,T]\\[4pt]
\si^-\,\gr n=\si^+\,\gr n,&&\mbox{on }\Gamma\times]0,T].
\end{eqnarray}
\end{subequations}
The solution to (\ref{eq:biot})
is then computed from the solution to the Green
Problem thanks to a  convolution by the source function. For instance we have:
$$P^+(x,y,t)=p^+(x,y,.)\ast f(.)=\int_0^t p^+(x,y,\tau)f(t-\tau)\,d\tau$$
(we have similar relations for the other
unknowns).  We also suppose that the poroelastic
medium is non dissipative, i.e the viscosity
$\eta^\pm=0$. Using the equations
(\ref{acousporo:eq:6},\ref{acousporo:eq:8}) we can
eliminate $\si^\pm$ and $p^\pm$ in
(\ref{biot_decomp}) and we obtain the equivalent
system:
\begin{equation}\label{syst-biot2}
\hspace*{-0.2cm}\left\{
\begin{array}{l}
\ro^\pm\,\ddu^\pm+\ro_f^\pm\,
\ddw^\pm-\alpha^\pm\,\nabla(\nabla\cdot \gr u_s^\pm)
+\mu^\pm\,\nabla\times(\nabla\times \gr u_s^\pm)-m^\pm\be^\pm\nabla(\nabla\cdot \gr w^\pm)\\[10pt]
=(f_u-\beta^+m^+f_p)\nabla(\delta_x\,\delta_y\,\delta_{z-h})\,\delta_t,\\[12pt]
\ro_f^\pm\,\ddus^\pm
+\ro_w^\pm\,\ddw^\pm-m^\pm\be^\pm\,\nabla(\nabla\cdot\gr u_s^\pm)
-m^\pm\,\nabla(\nabla\cdot\gr w^\pm)=(f_w-m^+f_p)\nabla(\delta_x\,\delta_y\,\delta_{z-h})\,\delta_t,
\end{array}
\right .
\end{equation}
with $\alpha^-=\la^-+2\mu^-+m^-{\be^-}^2$.\\\\
And the transmission conditions on $\Ga$  are rewritten as:
\begin{subequations}\label{eq:condi_trans_expli}
\begin{eqnarray}
&&u_{sx}^+=u_{sx}^-,\label{eq:condi_trans_expli1}\\[8pt]
&&u_{sy}^+=u_{sy}^-,\label{eq:condi_trans_expli2}\\[8pt]
&&u_{sz}^+=u_{sz}^-,\label{eq:condi_trans_expli3}\\[8pt]
&&\dsp w_{z}^-=w_z^+,\label{eq:condi_trans_expli4}\\[8pt]
&&m^+\be^+\,\nabla\cdot\gr u_s^++m^+\,\nabla\cdot\gr w^+=m^-\be^-\,\nabla\cdot\gr u_s^-+m^-\,\nabla\cdot\gr w^-,\label{eq:condi_trans_expli5}\\[12pt]
&&\dsp\mu^+(\p_zu_{sx}^++\p_x u_{sz}^+)=\mu^-(\p_zu_{sx}^-+\p_x u_{sz}^-),
\label{eq:condi_trans_expli6}\\[12pt]
&&\dsp\mu^+(\p_zu_{sy}^++\p_y u_{sz}^+)=\mu^-(\p_zu_{sy}^-+\p_y u_{sz}^-),
\label{eq:condi_trans_expli7}\\[12pt]
&&\dsp(\la^-+m^+{\be^+}^2)\nabla\cdot\gr u_s^++2\mu^+\p_z u_{sz}^++m^+\be^+\nabla\cdot\gr w^+=\label{eq:condi_trans_expli8}\\[6pt]
&&(\la^-+m^-{\be^-}^2)\nabla\cdot\gr u_s^-
 +2\mu^-\p_z u_{sz}^-+m^-\be^-\nabla\cdot\gr w^-.\nonumber
\end{eqnarray}
\end{subequations}
We split the displacement fields $\gr u_s^\pm$ and
$\gr w^\pm$ into irrotational and isovolumic fields
(P-wave and S-wave):
\begin{equation}\label{eq:irriso}
\gr u_s^\pm=\nabla \Theta_u^-+\nabla\times\Psi_u^\pm
\;\; ; \;\;\gr w^\pm= \nabla \Theta_w^\pm+\nabla\times\Psi_w^\pm.
\end{equation}
The vectors $\gr \Psi_u^\pm$ and $\gr \Psi_w^\pm$ are not uniquely 
defined since:
$$
\nabla\times (\gr \Psi_\ell^\pm+\nabla C)=\nabla\times\gr\Psi_\ell^\pm,\;\;\;\forall\,\ell\in\{u,w\}
$$  
for all scalar field $C$. To define a unique $\gr \Psi_\ell^\pm$ we impose the gauge condition:
$$
\nabla\cdot\gr \Psi_\ell^\pm=0.
$$  
The vectorial space of $\gr \Psi_\ell^\pm$ verifying this last condition is written as:
$$
\gr \Psi_\ell^\pm=\left[\begin{array}{c}
\p_y\\[5pt]
-\p_x\\[5pt]
0\end{array}\right]\Psi_{\ell,1}^\pm+
\left[\begin{array}{c}
\p_{xz}^2\\[5pt]
\p_{yz}^2\\[5pt]
-\p_{xx}^2-\p_{yy}^2\end{array}\right]\Psi_{\ell,2}^\pm,
$$ 
where $\Psi_{\ell,1}^\pm$ and $\Psi_{\ell,2}^\pm$ are two scalar fields. The displacement fields $\gr u_s^\pm$ and $\gr w^\pm$ are written in the form:
\begin{equation}\label{eq:changetp12}
\begin{array}{l}
\gr u_s^\pm=\nabla\Theta_u^\pm+\left[\begin{array}{c}
\p_{xz}^2\\[5pt]
\p_{yz}^2\\[5pt]
-\p_{xx}^2-\p_{yy}^2\end{array}\right]\Psi_{u,1}^\pm-\left[\begin{array}{c}
\p_y\\[5pt]
-\p_x\\[5pt]
0\end{array}\right]\Delta\Psi_{u,2}^\pm\\[35pt]
\gr w^\pm=\nabla\Theta_w^\pm+\left[\begin{array}{c}
\p_{xz}^2\\[5pt]
\p_{yz}^2\\[5pt]
-\p_{xx}^2-\p_{yy}^2\end{array}\right]\Psi_{w,1}^\pm-\left[\begin{array}{c}
\p_y\\[5pt]
-\p_x\\[5pt]
0\end{array}\right]\Delta\Psi_{w,2}^\pm.
\end{array}
\end{equation}
We can then rewrite system (\ref{syst-biot2}) in the following form:
\begin{equation}\label{equ:matrice}
\left\{
\begin{array}{ll}
A^+\ddot\Theta^+-B^+\Delta \Theta^+=
\delta_x\,\delta_y\,\delta_{z-h}\,\delta_t\,\gr F,&\mbox{in }\Om^+\times]0,T]\\[8pt]
A^-\ddot\Theta^--B^-\Delta \Theta^-=0,&\mbox{in }\Om^-\times]0,T]\\[8pt]
\ddot\Psi_{u,i}^\pm-{\Vs^\pm}^2\Delta \Psi_{u,i}^\pm=0,\ \ \i\in\{1,2\},&\mbox{in }\Om^\pm\times]0,T]\\[8pt]
\dsp
\ddot{\gr\Psi}_w^\pm=-\frac{\ro_f^\pm}{\ro_w^\pm}\ddot{\gr\Psi}_u^\pm,&\mbox{in }\Om^\pm\times]0,T]\end{array}
\right.
\end{equation}
where $\Theta^\pm=(\Theta_u^-,\Theta_w^-)^t$, 
$\gr F=(f_u-\beta^+m^+f_p,f_w-m^+f_p)^t$, $A^\pm$
and $B^\pm$ are $2\times 2$ symmetric matrices:
$$
A^\pm=\left(\begin{array}{cc}
\ro^\pm&\ro_f^\pm\\[8pt]
\ro_f^\pm&\ro_w^\pm
\end{array}\right)\;\;;\;\;B^\pm=\left(\begin{array}{cc}
\lambda^\pm+2\mu^\pm+m^\pm(\beta^\pm)^2&m^\pm\beta^\pm\\[8pt]
m^\pm\beta^\pm&m^\pm
\end{array}\right),
$$
and $$
\Vs^\pm=\sqrt{\frac{\mu\ro_w^\pm}{\ro^\pm\ro_w^\pm-{\ro_f^\pm}^2}}$$
is the S-wave velocity.\\\\
We multiply the first (resp. the second) equation of system (\ref{equ:matrice}) by
the inverse of $A^+$ (resp. $A^-$). The matrix  ${A^+}^{-1}B^+$ 
(resp. ${A^-}^{-1}B^-$) is diagonalizable:
${A^\pm}^{-1}B^\pm=\mat P^\pm D^\pm{\mat P^\pm}^{-1}$, where $\mat
P^\pm$ is the change-of-coordinates matrix,
$D^\pm=diag({\Vpf^\pm}^2,{\Vps^\pm}^2)$ is the diagonal matrix similar
to ${A^\pm}^{-1}B^\pm$, $\Vpf^\pm$ and $\Vps^\pm$ are
respectively the fast P-wave velocity and the slow
P-wave velocity ($\Vps^\pm<\Vpf^\pm$).\\\\
Using the change of variables: 
\begin{equation}\label{changethetaphi}
\Phi^\pm=(\Phi_{Pf}^\pm,\Phi_{Ps}^\pm)^t={\mat P^\pm}^{-1}\Theta^\pm,
\end{equation}
we obtain the uncoupled system on fast P-waves, slow P-waves and S-waves:
\begin{equation}\label{equ:matricediag}
\left\{
\begin{array}{ll}
\ddot\Phi^+-D^+\Delta \Phi^+=\delta_x\,\delta_y\,\delta_{z-h}\,\delta_t\,\gr F^+,&\mbox{in }\Om^+\times]0,T]\\[8pt]
\ddot\Phi^--D^-\Delta \Phi^-=0,&\mbox{in }\Om^-\times]0,T]\\[8pt]
\ddot\Psi_{u,\ell}^\pm-{\Vs^\pm}^2\Delta \Psi_{u,\ell}^\pm=0,\ \ \ell\in\{1,2\},&\mbox{in }\Om^\pm\times]0,T]\\[8pt]
\dsp
\gr\Psi_w^\pm=-\frac{\ro_f^\pm}{\ro_w^\pm}\gr\Psi_u^\pm,&\mbox{in }\Om^\pm\times]0,T]
\end{array}
\right.
\end{equation}
with $\gr F^+=(A^+\mat P^+)^{-1}\gr F=(F_{Pf}^+,F_{Ps}^+)^t$.\\\\
Using the transmission conditions  (\ref{eq:condi_trans_expli1})-(\ref{eq:condi_trans_expli2}), (\ref{eq:condi_trans_expli6})-(\ref{eq:condi_trans_expli7}) and the change of variables (\ref{eq:changetp12}), we obtain:
\begin{subequations}\label{eq:trans2eq}
\begin{eqnarray}
\p_x\Theta_u^++\p_{xz}^2\psi_{u,1}^+-\p_y(\Delta\psi_{u,2}^+)=
\p_x\Theta_u^-+\p_{xz}^2\psi_{u,1}^--\p_y(\Delta\psi_{u,2}^-),
&&\mbox{on }\Ga,\label{eq:trans2eq3}\\[18pt]
\p_y\Theta_u^++\p_{yz}^2\psi_{u,1}^++\p_x(\Delta\psi_{u,2}^+)=
\p_y\Theta_u^-+\p_{yz}^2\psi_{u,1}^-+\p_x(\Delta\psi_{u,2}^-),
&&\mbox{on }\Ga,\label{eq:trans2eq4}\\[18pt]
\begin{array}{l}
\mu^+\left(2\p_{xz}^2\Theta_u^++\p_x(\p_{zz}^2-\Delta_\perp)\Psi_{u,1}^+
-\p_{yz}^2\Delta\Psi_{u,2}^+\right)=\\[14pt]
\mu^-\left(2\p_{xz}^2\Theta_u^-+\p_x(\p_{zz}^2-\Delta_\perp)\Psi_{u,1}^-
-\p_{yz}^2\Delta\Psi_{u,2}^-\right),
\end{array}
&&\mbox{on }\Ga,\label{eq:trans2eq1}\\[18pt]
\begin{array}{l}
\mu^+\left(2\p_{yz}^2\Theta_u^++\p_y(\p_{zz}^2-\Delta_\perp)\Psi_{u,1}^+
+\p_{xz}^2\Delta\Psi_{u,2}^+\right)=\\[14pt]
\mu^-\left(2\p_{yz}^2\Theta_u^-+\p_y(\p_{zz}^2-\Delta_\perp)\Psi_{u,1}^-
+\p_{xz}^2\Delta\Psi_{u,2}^-\right)
\end{array}
&&\mbox{on }\Ga,\label{eq:trans2eq2}
\end{eqnarray}
\end{subequations}
with $\Delta_\perp=\p_{xx}^2+\p_{yy}^2$. Applying the  derivative $\p_y$ to the 
equation (\ref{eq:trans2eq3}) (resp. (\ref{eq:trans2eq1})), $\p_x$ to the equation (\ref{eq:trans2eq4})  (resp. (\ref{eq:trans2eq2})) and 
subtracting the first (resp. the third) obtained equation from the second (resp. the fourth) one, we get:
\begin{subequations}\label{eq:psi2}
\begin{eqnarray}
\Delta_\perp(\Delta\psi_{u,2}^+)=\Delta_\perp(\Delta\psi_{u,2}^-),&&\mbox{ on }\Ga\label{eq:psi22},\\[8pt]
\mu^+(\p_z\Delta_\perp)\Delta \Psi_{u,2}^+=\mu^-(\p_z\Delta_\perp)\Delta \Psi_{u,2}^-, &&\mbox{ on }\Ga,
\end{eqnarray}
\end{subequations}
moreover, using the third equation of (\ref{equ:matricediag}), we have $\Psi_{u,2}^\pm$ satisfies the wave equation:
$$
\ddot\Psi_{u,2}^\pm-{\Vs^\pm}^2\Delta \Psi_{u,2}^\pm=0, \;\;\mbox{in }\Om^\pm\times]0,T]
$$
and, since $\gr u_s^\pm$ and $\gr w^\pm$ satisfy, at  $t=0$, 
$
\gr u_s^\pm=\dot{\gr u}_s^\pm=\gr w^\pm=\dot{\gr w}^\pm=0,
$
we obtain: 
\begin{equation}\label{eq:psinull}
\Psi_{u,2}^\pm=0,\;\;\mbox{in }\Om^\pm\times]0,T], 
\end{equation}
and from (\ref{eq:trans2eq3})-(\ref{eq:trans2eq4}) we deduce the transmission condition equivalent to (\ref{eq:condi_trans_expli1}) and (\ref{eq:condi_trans_expli2}):
\begin{equation}\label{eq:transm_reduite1}
\p_x\Theta_u^++\p_{xz}^2\psi_{u,1}^+=
\p_x\Theta_u^-+\p_{xz}^2\psi_{u,1}^-,
\;\;\mbox{on }\Ga.
\end{equation}
In the same way, using the  equality (\ref{eq:psinull}), we can show 
that the two transmission conditions (\ref{eq:trans2eq1}) and (\ref{eq:trans2eq2}) 
are equivalent, which gives us: 
\begin{equation}\label{eq:transm_reduite2}
\mu^+\left(2\p_{xz}^2\Theta_u^++\p_x(\p_{zz}^2
-\Delta_\perp)\Psi_{u,1}^+\right)=
\mu^-\left(2\p_{xz}^2\Theta_u^-+\p_x(\p_{zz}^2-\Delta_\perp)\Psi_{u,1}^-\right),
\;\;\mbox{on }\Ga.
\end{equation}
We can then reduce the transmission conditions (\ref{eq:condi_trans_expli}) to 6 equations: (\ref{eq:transm_reduite1}, \ref{eq:condi_trans_expli3}, \ref{eq:condi_trans_expli4},\ref{eq:condi_trans_expli5},\ref{eq:transm_reduite2},\ref{eq:condi_trans_expli8}).\\\\
Finally, we obtain the Green problem equivalent to (\ref{biot_decomp}):
\begin{equation}\label{equ:diagsyst}
\left\{
\begin{array}{ll}
\ddot \Phi_i^+-{V_i^+}^2\Delta\Phi_i^+=\delta_x\,\delta_y\,\delta_{z-h}\,\delta_t\,F_i^+,\quad i\in\{Pf,Ps\}&z>0\\[8pt]
\ddot \Phi_{S}^+-{V_{S}^+}^2\Delta\Phi_S^+=0&z>0\\[8pt]
\ddot\Phi_i^--{V_i^-}^2\Delta \Phi_i^-=0,\quad i\in\{Pf,Ps,S\}&z<0\\[8pt]
\dsp {\cal B} (\Phi_{Pf}^+,\Phi_{Ps}^+,\Phi_S^+,\Phi_{Pf}^-,\Phi_{Ps}^-,\Phi_S^-)=0,&z=0
\end{array}
\right.
\end{equation}
where we have set $\Phi_S^\pm=\Psi_{u,1}^\pm$ in order to have similar notations for the $Pf$, $Ps$ and $S$ waves.  The operator ${\cal B}$ represents the
transmission conditions on $\Gamma$:
$$
{\cal B} \left(
\begin{array}{l}
\Phi_{Pf}^+\\[4pt]
\Phi_{Ps}^+\\[4pt]
\Phi_S^+\\[4pt]
\Phi_{Pf}^-\\[4pt]
\Phi_{Ps}^-\\[4pt]
\Phi_S^-
\end{array}\right)=\left[
\begin{array}{cccccc}
\mat P_{11}^+\,\p_x&\mat P_{12}^+\,\p_x&\p_{xz}^2&-\mat P_{11}^-\,\p_x&
-\mat P_{12}^-\,\p_x&-\p_{xz}^2\\[8pt]
\mat P_{11}^+\,\p_z&\mat P_{12}^+\,\p_z&-\Delta_\perp&-\mat P_{11}^-\,\p_z&
-\mat P_{12}^-\,\p_z&\Delta_\perp\\[8pt]
\mat P_{21}^+\,\p_z&\mat P_{22}^+\,\p_z&\dsp\frac{\ro_f^+}{\ro_w^+}\Delta_\perp&-\mat P_{21}^-\,\p_z&
-\mat P_{22}^-\,\p_z&\dsp-\frac{\ro_f^-}{\ro_w^-}\Delta_\perp\\[8pt]
{\cal B}_{41}&{\cal B}_{42}&0&{\cal B}_{44}&{\cal B}_{45}&0\\[8pt]
{\cal B}_{51}&{\cal B}_{52}&\mat B_{53}
&{\cal B}_{54}&{\cal B}_{55}&\mat B_{56}\\[8pt]
{\cal B}_{61}&{\cal B}_{62}&-2\mu^+\p_{z}(\Delta_\perp)&{\cal B}_{64}&{\cal B}_{65}&
2\mu^-\p_z(\Delta_\perp)
\end{array}
\right]\left[
\begin{array}{l}
\Phi_{Pf}^+\\[9pt]
\Phi_{Ps}^+\\[9pt]
\Phi_S^+\\[9pt]
\Phi_{Pf}^-\\[9pt]
\Phi_{Ps}^-\\[9pt]
\Phi_S^-
\end{array}
\right]
$$
where $\mat P_{ij}^\pm$, $i,j=1,2$ are the components of the 
change-of-coordinates matrix $\mat P^\pm$ and 
$$
\begin{array}{l}
\dsp{\cal B}_{41}=\frac{m^+(\be^+\mat P_{11}^+
+\mat P_{21}^+)}{{\Vpf^+}^2}\,\p_{tt}^2\ ;\ 
\dsp{\cal B}_{42}=\frac{m^+(\be^+\mat P_{12}^++\mat P_{22}^+)}{{\Vps^+}^2}
\,\p_{tt}^2;\\[18pt]
\dsp{\cal B}_{44}=-\frac{m^-(\be^-\mat P_{11}^
-+\mat P_{21}^-)}{{\Vpf^-}^2}\,\p_{tt}^2\ ;\ 
\dsp{\cal B}_{45}=-\frac{m^-(\be^-\mat P_{12}^-
+\mat P_{22}^-)}{{\Vps^-}^2}\,\p_{tt}^2;\\[18pt]
{\cal B}_{51}=2\mu^+\mat P_{11}^+\,\p_{xz}^2\ ; \ {\cal B}_{52}=2\mu^+\mat P_{12}^+\,\p_{xz}^2\ ;\ \mat B_{53}=\mu^+\p_x(\p_{zz}^2-\Delta_\perp);\\[12pt]
{\cal B}_{54}=-2\mu^-\mat P_{11}^-\,\p_{xz}^2\ ;
\ {\cal B}_{55}=-2\mu^-\mat P_{12}^-\,\p_{xz}^2\ ; \ 
\mat B_{56}=-\mu^-\p_x(\p_{zz}^2-\Delta_\perp);\\[12pt]
\dsp {\cal B}_{61}=\frac{(\la^++m^+{\be^+}^2)\mat P_{11}^++
m^+\be^+\mat P_{21}^+}{{\Vpf^+}^2}\p_{tt}^2
+2\mu^+\mat P_{11}^+\p_{zz}^2;\\[12pt] 
\dsp {\cal B}_{62}=\dsp \frac{(\la^++m^+{\be^+}^2)\mat P_{12}^++
m^+\be^+\mat P_{22}^+}{{\Vps^+}^2}\,\p_{tt}^2
+2\mu^+\mat P_{12}^+\p_{zz}^2;\\[12pt]
\dsp {\cal B}_{64}=\frac{(\la^-+m^-{\be^-}^2)\mat P_{11}^-+
m^-\be^-\mat P_{21}^-}{{\Vpf^-}^2}\p_{tt}^2
+2\mu^-\mat P_{11}^-\p_{zz}^2;\\[12pt] 
\dsp {\cal B}_{65}=\dsp \frac{(\la^-+m^-{\be^-}^2)\mat P_{12}^-+
m^-\be^-\mat P_{22}^-}{{\Vps^-}^2}\,\p_{tt}^2+2\mu^-\mat P_{12}^-\p_{zz}^2.
\end{array}
$$
To obtain this operator we have used the reduced transmission conditions (\ref{eq:transm_reduite1}, \ref{eq:condi_trans_expli3}, \ref{eq:condi_trans_expli4},\ref{eq:condi_trans_expli5},\ref{eq:transm_reduite2},\ref{eq:condi_trans_expli8}), the change of variables (\ref{eq:irriso},\ref{changethetaphi}) and the uncoupled system (\ref{equ:matricediag}).\\\\
Moreover, from the unknowns $\Phi_{Pf}^\pm$, $\Phi_{Ps}^\pm$ and $\Phi_S^\pm$ we can determine the solid displacement $\gr u_s^\pm$ and the  relative displacement $\gr w^\pm$ by using the change of variables presented below.
\section{Expression of the analytical solution}
Since the problem is invariant by a rotation
around the $z$-axis, we will only consider the
case $y=0$ and $x>0$, so that the y-component of all the
displacements are zero. The solution for $y\neq0$
or $x\leq0$
is deduced from the solution for $y=0$ by the
relations
\begin{eqnarray}
  \label{acousporo3d:eq:2}
&&\dsp
u^\pm_{s\,x}(x,y,z,t)=\frac{x}{\sqrt{x^2+y^2}}u^\pm_{s\,x}(\sqrt{x^2+y^2},0,z,t)\\[10pt]
&&\dsp
u^\pm_{s\,y}(x,y,z,t)=\frac{y}{\sqrt{x^2+y^2}}u^\pm_{s\,x}(\sqrt{x^2+y^2},0,z,t)\\[10pt]
&&\dsp u^\pm_{s\,z}(x,y,z,t)=u^\pm_{s\,z}(\sqrt{x^2+y^2},0,z,t)
\end{eqnarray}
To state our results, we need the following notations and definitions:
\begin{enumerate}
\item {\bf Definition of the complex square
    root}. 
 For $q\in\setC\backslash\setR^-$, we use the following definition of the square root $g(q)=q^{1/2}$:
\[
  g(q)^2=q \quad \hbox{ and } \quad
  \Re e[g(q)]>0.
\]
The branch cut of $g(q)$ in the complex plane will thus be the half-line defined by $\{q \in \setR^- \}$ (see Fig.~\ref{simple:fig:15}).
In the following, we use the abuse of notation $g(q)=\ic\sqrt{-q}$ for $q\in\setR^-$.
\\[10pt]
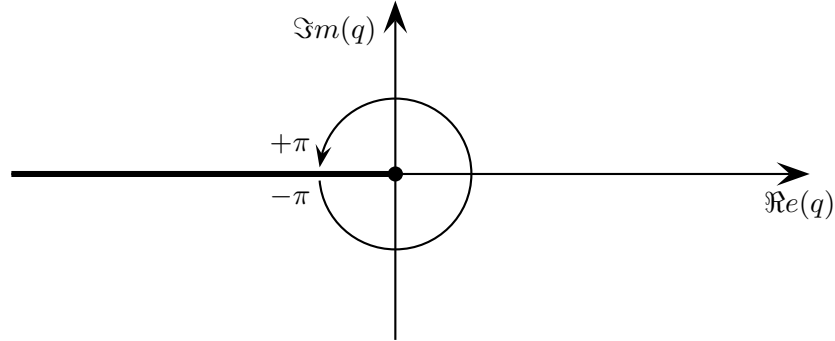
\begin{figure}[htbp]
\setlength{\unitlength}{1cm}
\psset{xunit=1cm,yunit=1cm,runit=1cm}
\begin{picture}(10,5)
\put(0,-2){\put(6.7,6.){$\Im m(q)$}
\put(12.9,3.7){$\Re e(q)$}
\put(6.4,3.8){$-\pi$}
\put(6.4,4.5){$+\pi$}
\psline[arrows=->,arrowsize=.2 4](3,4.2)(13.5,4.2)
\psline[arrows=->,arrowsize=.2 4](8.05,2)(8.05,6.5)
\psline[linewidth=0.08](3,4.2)(8.05,4.2)
\pscircle[fillstyle=solid,fillcolor=black](8.05,4.2){.1}
\psarc[arrows=->,arrowsize=.1 4](8.05,4.2){1}{-175}{175}}
\end{picture}
\caption{Definition of the function $x\mapsto (x)^{1/2}$}
\label{simple:fig:15}
\end{figure}
\item {\bf Definition of the fictitious velocities}
For a given $q\in\setR$, we define the fictitious
velocities $\mat V_i^\pm(q)$ for
$i\in\{Pf,Ps,S\}$ by 
$$\mat V^\pm_i:=\mat V^\pm_i(q)=V^\pm_i\sqrt{\frac{1}{1+{V^\pm_i}^2q^2}}.$$
These fictitious velocities will be helpful to turn
the 3D-problem into the sum of 2D-problems indexed
by the variable $q$. Note that
$\mat V^\pm_i(0)$ correspond to the real
velocities $V^\pm_i$.
\item {\bf Definition of the functions $\kappa^{\pm}_{i}$}.
For $i\in\{Pf,Ps,S\}$
  and $(q_x,q_y)\in\setC\times\setR$, we define the functions
$$
\kappa_i^\pm:=\kappa_i^\pm(q_x,q_y)=\left(\frac{1}{{V_i^\pm}^2}+q_x^2+q_y^2\right)^{1/2}=\left(\frac{1}{{\mat V_i^\pm}^2(q_y)}+q_x^2\right)^{1/2}.$$
\item {\bf Definition of the reflection and transmission coefficients.}
For a given $\gr q =(q_x,q_y)\in\setC\times\setR$, we denote by $\rff(\gr q)$, $\rfs(\gr q)$,  $\rfpsi(\gr q)$, 
 $\tff(\gr q)$, $\tfs(\gr q)$, and 
$\tfpsi(\gr q)$ 
the solution to the linear system
\begin{equation*}
\hspace*{-.5cm}\mat {\cal A}(\gr q)
\left[
  \begin{array}{l}
    \rff(\gr q) \\[10pt]
    \rfs(\gr q) \\[10pt]
    \rfpsi(\gr q) \\[10pt]
 \tff(\gr q)\\[10pt]
 \tfs(\gr q)\\[10pt]
\tfpsi(\gr q)
  \end{array}
\right]
=\frac{1}{2\kafp(\gr q){V^+_{Pf}}^2}
\left[
  \begin{array}{c}
\dsp     \ic q_x\mat P_{11}^+\\[10pt]
\dsp  -\kafp(\gr q) \mat P_{11}^+\\[10pt]
\dsp -\kafp(\gr q) \mat P_{21}^+\\[10pt]
\dsp -\frac{m^+}{{V^+_{Pf}}^2}(\beta^+\mat P_{11}^++\mat P_{21}^+)\\[15pt]
2\ic q_x\mu^+\kafp(\gr q)\mat P_{11}^+\\[10pt]
\dsp -\frac{(\la^++m^+{\be^+}^2)\mat P_{11}^++
m^+\be^+\mat P_{21}^+}{{\Vpf^+}^2}-2\mu^+\mat P_{11}^+\kafp^2(\gr q)
  \end{array}
\right]
\end{equation*}
 and by $\rsf(\gr q)$, $\rss(\gr q)$,  $\rspsi(\gr q)$, $\tsf(\gr q)$, $\tss(\gr q)$ and  $\tspsi(\gr q)$ the solution to the linear system
\begin{equation*}
\hspace*{-.5cm}\mat {\cal A}(\gr q)
\left[
  \begin{array}{l}
    \rsf(\gr q) \\[10pt]
    \rss(\gr q) \\[10pt]
    \rspsi(\gr q) \\[10pt]
 \tsf(\gr q)\\[10pt]
 \tss(\gr q)\\[10pt]
\tspsi(\gr q)
  \end{array}
\right]
=\frac{1}{2\kasp(\gr q){V^+_{Ps}}^2}
\left[
  \begin{array}{c}
\dsp     \ic q_x\mat P_{12}^+\\[10pt]
\dsp  -\kasp(\gr q) \mat P_{12}^+\\[10pt]
\dsp -\kasp(\gr q) \mat P_{22}^+\\[10pt]
\dsp -\frac{m^+}{{V^+_{Ps}}^2}(\beta^+\mat P_{12}^++\mat P_{22}^+)\\[15pt]
2\ic q_x\mu^+\kasp(\gr q)\mat P_{12}^+\\[10pt]
\dsp -\frac{(\la^++m^+{\be^+}^2)\mat P_{12}^++
m^+\be^+\mat P_{22}^+}{{\Vps^+}^2}-2\mu^+\mat P_{12}^+\kasp^2(\gr q)
  \end{array}
\right]
,
\end{equation*}
where the matrix $\mat A(\gr q)$ is defined
for $\gr q =(q_x,q_y)\in\setC\times\setR$ by:
$$
\hspace*{-.9cm} A(\gr q)=\left[\begin{array}{cccccc}
 -\ic q_x\mat P^+_{11}&-\ic q_x \mat P^+_{12}&\ic q_x\kapsip(\gr q)&\ic q_x\mat P^-_{11}&\ic q_x \mat P^-_{12}&\ic q_x\kapsi(\gr q)\\[10pt]
-\kafp(\gr q) \mat P^+_{11}&-\kasp (\gr q)\mat P^+_{12}&q_x^2+q_y^2 &-\kaf (\gr q)\mat P^-_{11}&-\kas (\gr q)\mat P^-_{12}&-q_x^2-q_y^2\\[10pt]
-\kafp(\gr q) \mat P^+_{21}&-\kasp (\gr q)\mat P^+_{22}&-(q_x^2+q_y^2)\frac{\rho_f^+}{\rho_w^+}&-\kaf(\gr q) \mat P^-_{21}&-\kas(\gr q) P^-_{22}&(q_x^2+q_y^2)\frac{\rho_f^-}{\rho_w^-}\\[10pt]
 \mat A_{41}(\gr q)&\mat A_{42}(\gr q)&0&\mat A_{44}(\gr q)&\mat A_{45}(\gr q)&0\\[10pt]
\mat A_{51}(\gr q)&\mat A_{52}(\gr q)&\mat A_{53}(\gr q)&\mat A_{54}(\gr q)&\mat A_{55}(\gr q)&\mat A_{56}(\gr q)\\[10pt]
\mat A_{61}(\gr q)&\mat A_{62}(\gr q)& A_{63}(\gr q) &\mat A_{64}(\gr q)&\mat A_{65}(\gr q)& \mat A_{66}(\gr q)
\end{array}\right],
$$
with
$$\begin{array}{l}
\begin{array}{rclrcl}
\dsp \mat A_{41}(\gr q)&=&\dsp \frac{m^+}{{V_{Pf}^+}^2}\left[\beta^+\mat P^+_{11}+\mat P^+_{21}\right];&\dsp \mat A_{44}(\gr q)&=&\dsp -\frac{m^-}{{V_{Pf}^-}^2}\left[\beta^-\mat P^-_{11}+\mat P^-_{21}\right];
\\[18pt]
\dsp \mat A_{42}(\gr q)&=&\dsp \frac{m^+}{{V_{Ps}^+}^2}\left[\beta^+\mat P^+_{12}+\mat P^+_{22}\right];&\dsp \mat A_{45}(\gr q)&=&\dsp -\frac{m^-}{{V_{Ps}^-}^2}\left[\beta^-\mat P^-_{12}+P^-_{22}\right];
\end{array}
\\[35pt]
\begin{array}{rclrcl}
\dsp \mat A_{51}(\gr q)&=&\dsp 2\ic q_x\mu^+\kafp(\gr q)  \mat P^+_{11};&
\dsp A_{54}(\gr q)&=&2\ic q_x\mu^-\kaf(\gr q) \mat P^-_{11};
\\[10pt]
\dsp \mat A_{52}(\gr q)&=&\dsp 2\ic q_x\mu^+\kasp(\gr q)  \mat P^+_{12};&
\dsp \mat A_{55}(\gr q)&=&\dsp 2\ic q_x \mu^-\kas(\gr q)  P^-_{12};\\[10pt]
\dsp  \mat A_{53}(\gr q)&=&\dsp -\ic \mu^+q_x({\kapsip}^2(\gr q)+q_x^2+q_y^2);
&\dsp \mat A_{56}(\gr q)&=&\dsp \ic \mu^-q_x({\kapsi}^2(\gr q)+q_x^2+q_y^2);
\end{array}
\\[35pt]
\begin{array}{l}
\dsp \mat A_{61}(\gr q)=\frac{(\lambda^++m^+{\beta^+}^2)\mat P^+_{11}+m^+\beta^+\mat P^+_{21}}{{V_{Pf}^+}^2}+2\mu^+\kafp^2(\gr q)\mat P^+_{11};\\[18pt]
\dsp \mat A_{62}(\gr q)=\frac{(\lambda^++m^+{\beta^+}^2)\mat P^+_{12}+m^+\beta^+\mat P^+_{22}}{{V_{Ps}^+}^2}+2\mu^+{\kasp(\gr q)}^2\mat P^+_{12};\\[18pt]
\dsp \mat A_{63}(\gr q)=-2(q_x^2+q_y^2)\mu^+\kapsip(\gr q);\\[18pt]
\dsp \mat A_{64}(\gr q)=-\frac{(\lambda^-+m^-{\beta^-}^2)\mat P^-_{11}+m^-\beta^- P^-_{21}}{{V_{Pf}^-}^2}-2\mu^-{\kaf}^2(\gr q)\mat P^-_{11};\\[18pt]
\dsp \mat A_{65}(\gr q)=-\frac{(\lambda^-+m^-{\beta^-}^2)\mat P^-_{12}+m^-\beta^- P^-_{22}}{{V_{Ps}^-}^2}-2\mu^-{\kas}^2(\gr q)\mat P^-_{12};\\[18pt]
\dsp  \mat A_{66}(\gr q)=-2(q_x^2+q_y^2)\mu^-\kapsi (\gr q).
\end{array}
\end{array}
$$
\end{enumerate}
We also denote by $V_{\max}$ the greatest velocity in
the medium: $$V_{\max}=\max(\Vpf^+,\Vps^+,\Vs^+,\Vpf^-,\Vps^-,\Vs^-).$$

We can now present the expression of the solution
to the Green Problem:
\begin{theo}
The solid displacement in the top medium
is given by
$$
\begin{array}{rcccccccc}
u_s^+(x,0,z,t)&=&\gr u^+_{Pf}(x,z,t)
 &+&\gr u^+_{PfPf}(x,z,t)&+&\gr u^+_{PfPs}(x,z,t)&+&\gr u^+_{PfS}(x,z,t)\\[10pt]
 &+&\gr u^+_{Ps}(x,z,t)&+&\gr u^+_{PsPf}(x,z,t)&+&\gr u^+_{PsPs}(x,z,t)
 &+&\gr u^+_{PsS}(x,z,t)
  \end{array}
$$
and the 
solid displacement in the bottom medium is given by
$$\begin{array}{rcccccc}
\gr  u_s^-(x,0,z,t)&=& \gr u^+_{PfPf}(x,z,t)&+&\gr u^+_{PfPs}(x,z,t)&+&\gr u^+_{PfS}(x,z,t)\\[10pt]
&+&\gr u^+_{PsPf}(x,z,t)&+&\gr u^+_{PsPs}(x,z,t)&+&\gr u^+_{PsS}(x,z,t),
\end{array}
$$ where
  \begin{itemize}
  \item $\gr u^+_{Pf}$ is the solid displacement of the
    incident $Pf$ wave and satisfies:
$$\left\{
\begin{array}{lcl}
u^+_{Pf,x}(x,z,t)&:=&-\dsp\frac{\mat P^+_{11}F^+_{Pf}}{ {V_{Pf}^+}^2}
\frac{xtH(t-t_0)}{4\pi r^3}\\[20pt]
\dsp  u^+_{Pf,z}(x,z,t)&:=&-\dsp\frac{\mat P^+_{11}F^+_{Pf}}{ {V_{Pf}^+}^2}
\frac{(z-h)tH(t-t_0)}{4\pi r^3}\
\end{array}\right.,$$
where $H$ denotes the usual Heaviside function. Moreover we set
$r=(x^2+(z-h)^2)^{1/2}$ and $t_0=r/V_{Pf}^+$ denotes the time
arrival of the incident $Pf$ wave at point $(x,0,z)$.
\\[5pt]
  \item $\gr u^+_{Ps}$ is the solid displacement of the
    incident $Ps$ wave and satisfies:
$$\left\{
\begin{array}{lcl}
u^+_{Ps,x}(x,z,t)&:=&-\dsp\frac{\mat P^+_{12}F^+_{Ps}}{ {V_{Ps}^+}^2}
\frac{xtH(t-t_0)}{4\pi r^3}\\[20pt]
\dsp  u^+_{Ps,z}(x,z,t)&:=&-\dsp\frac{\mat P^+_{12}F^+_{Ps}}{ {V_{Ps}^+}^2}
\frac{(z-h)tH(t-t_0)}{4\pi r^3}.
\end{array}\right.$$
We set here
$r=(x^2+(z-h)^2)^{1/2}$ and $t_0=r/V_{Ps}^+$ denotes the time
arrival of the incident $Ps$ wave at point $(x,0,z)$.
\\[5pt]
  \item  $\gr u^+_{PfPf}$ is the solid displacement of the
    reflected $PfPf$ wave (the $Pf$ reflected wave generated by the 
$Pf$ incident wave) and satisfies:
$$\left\{
\begin{array}{lcl}
 \dsp  u_{PfPf,x}^+(x,z,t)&=&\dsp\mat P^+_{11}F^+_{Pf}\int_{0}^{q_1(t)} \frac{\Im
   m\Big[\ic\upsilon(t,q)\kafp(\upsilon(t,q))\rff(\upsilon(t,q))\Big]}{\pi^2r \sqrt{q^2+q_0^2(t)}}dq,
\\[15pt]
\dsp u_{PfPf,z}^+(x,z,t)&=&\dsp \mat P^+_{11}F^+_{Pf}
\int_{0}^{q_1(t)}\frac{\Im
  m\Big[\kafp^2(\upsilon(t,q))\rff(\upsilon(t,q))\Big]}{\pi^2r \sqrt{q^2+q_0^2(t)}}dq,
\end{array}\right.$$
$\dsp\hbox{if }t_{\hbox{h}_1} <t\leq t_{0} \hbox{ and } \frac{x}{r}>\frac{V_{Pf}^+}{V_{\max}}$,
$$\left\{
\begin{array}{lcl}
  \dsp u_{PfPf,x}^+(x,z,t)&=&\dsp 
 \begin{array}[t]{ll}
&\dsp \mat P^+_{11}F^+_{Pf} \int^{q_1(t)}_{q_0(t)}\frac{\Im
   m\Big[(\ic\upsilon(t,q)\kafp(\upsilon(t,q))\rff(\upsilon(t,q))\Big]}{\pi^2r \sqrt{q^2-q_0^2(t)}}dq
\\[18pt]
 -&\dsp\mat P^+_{11}F^+_{Pf}\int_{0}^{q_0(t)}\frac{\Re
   e\Big[\ic\gamma(t,q)\kafp(\gamma(t,q))\rff(\gamma(t,q))\Big]}{\pi^2r \sqrt{q_0^2(t)-q^2}}dq,
\end{array}
\\[60pt]
\dsp u_{PfPf,z}^+(x,z,t)&=&\dsp 
  \begin{array}[t]{ll}
&\dsp\mat P^+_{11}F^+_{Pf}\int^{q_1(t)}_{q_0(t)}\frac{\Im
  m\Big[\kafp^2(\upsilon(t,q))\rff(\upsilon(t,q))\Big]}{\pi^2r \sqrt{q^2-q_0^2(t)}}dq
\\[15pt]
-&\dsp \mat P^+_{11}F^+_{Pf}\int_{0}^{q_0(t)} \frac{\Re
  e\Big[\kafp^2(\gamma(t,q))\rff(\gamma(t,q))\Big]}{\pi^2r \sqrt{q_0^2(t)-q^2}}dq,
\end{array}
\end{array}\right.$$
$\dsp\hbox{if }t_{0} <t\leq t_{\hbox{h}_2}  \hbox{ and } \frac{x}{r}>\frac{V_{Pf}^+}{V_{\max}}$,
$$\left\{
\begin{array}{lcl}
  \dsp u_{PfPf,x}^+(x,z,t)&=&-\mat P^+_{11}F^+_{Pf}\dsp  \int^{q_0(t)}_{0}\frac{\Re e\Big[\ic\gamma(t,q)\kafp(\gamma(t,q))\rff(\gamma(t,q))\Big]}{\pi^2r \sqrt{q_0^2(t)-q^2}}dq,
\\[15pt]
\dsp u_{PfPf,z}^+(x,z,t)&=&-\mat P^+_{11}F^+_{Pf}\dsp \int^{q_0(t)}_{0} \frac{\Re e\Big[\kafp^2(\gamma(t,q))\rff(\gamma(t,q))\Big]}{\pi^2r \sqrt{q_0^2(t)-q^2}}dq,
\end{array}\right.$$
$\dsp\hbox{if } t_{\hbox{h}_2}<t
\hbox{ and } \frac{x}{r}>\frac{V_{Pf}^+}{V_{\max}}$ or
 $\dsp\hbox{if }t_{0} <t
\hbox{ and } \frac{x}{r}\leq \frac{V_{Pf}^+}{V_{\max}}$
and 
$$
 \gr
u_{PfPf}(x,z,t)=0 \hbox{ else }.$$
We set here $r=(x^2+(z+h)^2)^{1/2}$  and
$t_0=r/V_{Pf}^+$ denotes the arrival time of the reflected $PfPf$
volume wave at point $(x,0,z)$, 
\begin{eqnarray}
  t_{h_1}&=&(z+h)\sqrt{\frac{1}{{V_{Pf}^+}^2}-\frac{1}{V^2_{\max}}}+\frac{|x|}{V_{\max}}
\end{eqnarray}
denotes the arrival time of the reflected $PfPf$
head-wave at point $(x,0,z)$ and 
\begin{eqnarray}
  t_{h_2}&=&\frac{r}{z+h}\sqrt{\frac{1}{{V_{Pf}^+}^2}-\frac{1}{{V^2_{\max}}}}
\end{eqnarray}
denotes the time after which there is no longer head
wave at point $(x,0,z)$, (contrary to the 2D case, this time does not
coincide with the arrival time of the volume wave).
We also define the functions $\gamma$, $\upsilon$,
 $q_0$ and $q_1$ by
$$\gamma : \{t\in\setR \,|\, t>t_0\}\times\setR\mapsto
\setC:=
\gamma(t,q_y)=\ic \frac
{xt}{r^2}+\frac{z+h}{r}\sqrt{\frac{t^2}{r^2}-\frac{1}{{\mat
      V_{Pf}^{+2}}(q_y)}}$$
$$\upsilon : \{t\in\setR \,|\, t_{h_1}<t<t_{h_2}\}\times\setR\mapsto
\setC:=
\upsilon(t,q_y)=-\ic\left(\frac{z+h}{r}-\sqrt{\frac{1}{\mat
      V_{Pf}^{+2}(q_y)}-\frac{t^2}{r^2}}+\frac
  {x}{r^2}t\right),$$
$$q_0:\setR\to\setR:=q_0(t)=\sqrt{\left|\frac{t^2}{r^2}-\frac{1}{
      {V_{Pf}^+}^2}\right|}$$
and 
$$q_1:\setR\to\setR:=q_1(t)=\sqrt{\frac{1}{x^2}\left(t-(z+h)\sqrt{\frac{1}{{V_{Pf}^+}^2}-\frac{1}{V^2_{\max}}}\right)^2-\frac{1}{V_{\max}^2}}.$$
  \item  $\gr u^+_{PfPs}$ is the solid displacement of the
    reflected $PfPs$ wave and satisfies:
$$
\left\{
\begin{array}{ll}
 \dsp u^+_{PfPs,x}(x,z,t)=-\frac{\mat
   P_{12}^+F^+_{Pf}}{\pi^2}\int_{0}^{q_{1}(t)}
\Re
 e\left[\ic\upsilon(t,q)\rfs(\upsilon(t,q))\frac{\partial \upsilon}{\partial t}(t,q)\right]dq,
\\[18pt]
\dsp u^+_{PfPs,z}(x,z,t)= 
-\frac{\mat P_{12}^+F^+_{Pf}}{\pi^2}\int_{0}^{q_{1}(t)}
\Re
e\left[\kasp(\upsilon(t,q))\rfs(\upsilon(t,q))\frac{\partial\upsilon}{\partial    t}(t,q)\right]dq,
\end{array}\right.$$
$\dsp\hbox{if }t_{h_1} <t\leq t_{0} \hbox{ and } \left|\Im
    m\left[\gamma(t_{0},0)\right]\right|
  <\frac{1}{V_{\max}},$
$$\left\{
\begin{array}{lcl}
  \dsp u^+_{PfPs,x}(x,z,t)&=&\dsp 
\begin{array}[t]{ll}
-&\dsp\frac{\mat
    P_{12}^+F^+_{Pf}}{\pi^2}
\int_{0}^{q_0(t)}
\Re
e\left[\ic\gamma(t,q)\rfs(\gamma(t,q))\frac{\partial
    \gamma}{\partial t}(t,q)\right]dq
\\[18pt]
-&\dsp\frac{\mat
    P_{12}^+F^+_{Pf}}{\pi^2}
\int_{q_0(t)}^{q_{1}(t)}
\Re
 e\left[\ic\upsilon(t,q)\rfs(\upsilon(t,q))\frac{\partial \upsilon}{\partial t}(t,q)\right]dq,
\end{array}
\\[60pt]
\dsp u^+_{PfPs,z}(x,z,t)&=&\dsp 
\begin{array}[t]{ll}
-&\dsp \frac{\mat
  P_{12}^+F^+_{Pf}}{\pi^2}\int_{0}^{q_0(t)}
\Re
e\left[\kasp(\gamma(t,q))\rfs(\gamma(t,q))\frac{\partial
    \gamma
}{\partial t}(t,q)\right]dq
\\[18pt]
-&\dsp\frac{\mat
  P_{12}^+F^+_{Pf}}{\pi^2}
\int_{q_0(t)}^{q_{1}(t)}
\Re
e\left[\kasp(\upsilon(t,q))\rfs(\upsilon(t,q))\frac{\partial\upsilon}{\partial    t}(t,q)\right]dq,
\end{array}
\end{array}\right.$$
$\dsp\hbox{if }t_{0} <t\leq t_{h_2}  \hbox{ and }  \left|\Im
    m\left[\gamma(t_{0},0)\right]\right|
  <\frac{1}{V_{\max}}$,
$$
\left\{
\begin{array}{ll}
  \dsp u^+_{PfPs,x}(x,z,t)=- \frac{\mat P_{12}^+F^+_{Pf}}{\pi^2}\int_{0}^{q_0(t)}
\Re
e\left[\ic\gamma(t,q)\rfs(\gamma(t,q))\frac{\partial
    \gamma}{\partial t}(t,q)\right]dq,
\\[18pt]
\dsp u^+_{PfPs,z}(x,z,t)=-\frac{\mat
  P_{12}^+F^+_{Pf}}{\pi^2}\int_{0}^{q_0(t)}
\Re
e\left[\kasp(\gamma(t,q))\rfs(\gamma(t,q))\frac{\partial
    \gamma
}{\partial t}(t,q)\right]dq,
\end{array}\right.$$
$\dsp\hbox{if } t_{h_2}<t \hbox{ and } \left|\Im
  m\left[\gamma(t_{0},0)\right]\right|
<\frac{1}{V_{\max}}$ or $\dsp\hbox{if }t_{0} <t
\hbox{ and } \left|\Im
  m\left[\gamma(t_{0},0)\right]\right|
\geq\frac{1}{V_{\max}}$\\
and $\gr
u^+_{PfPs}(x,z,t)=0$ else.\\\\
$t_0$ denotes here the arrival time of the
reflected $PfPs$ volume wave at point
$(x,0,z)$ (its calculation is similar to the
calculation of the arrival time of the transmitted
wave, see the appendix of~\cite{RAP2}),
\begin{eqnarray}
  t_{h_1}&=&h\sqrt{\frac{1}{{V^+_{Pf}}^2}-\frac{1}{V^2_{\max}}}+z\sqrt{\frac{1}{{V^+_{Ps}}^2}-\frac{1}{V^2_{\max}}}+\frac{|x|}{V_{\max}}
\end{eqnarray}
denotes the arrival time of the reflected $PfPs$
head wave at point $(x,0,z)$,
\begin{eqnarray}
  t_{h_2}&=&\frac{h^2+z^2+hz\left(\frac{c_2}{c_1}+
\frac{c_1}{c_2}
\right)+x^2}
{\frac{h}{c_1}+\frac{z}{c_2}}
\end{eqnarray}
 denotes the time after which there is no longer head
wave at point $(x,0,z)$, where
$$c_1=\sqrt{\frac{1}{{V_{Pf}^+}^2}-\frac{1}{{V^2_{\max}}}} \hbox{ and } c_2=\sqrt{\frac{1}{{V^+_{Ps}}^2}-\frac{1}{{V^2_{\max}}}}.$$
The function $q_0:[t_0\,;\,+\infty]\mapsto\setR^+$ is
the reciprocal function of $\tilde t_0
:\setR^+\mapsto:[t_0,+\infty]$, where $\widetilde
t_0(q)$ is the arrival time  at point $(x,0,z)$ of the fictitious 
transmitted $PfPs$ volume wave, propagating at a velocity $\mat
V_{Pf}^+(q)$ from the source to the interface and at velocity $\mat
V^+_{Ps}(q)$ from the interface to point $(x,0,z)$ (we refer again to~\cite{RAP2} for details on its calculation).\\[5pt]
The function $q_1:[t_1\,;\,t_0]\mapsto\setR^+$ is defined by 
$$q_{1}(t)=\sqrt{\frac{1}{x^2}\left(t-z\sqrt{\frac{1}{{V^+_{Ps}}^2}-\frac{1}{V^2_{\max}}}
-h\sqrt{\frac{1}{{V_{Pf}^+}^2}-\frac{1}{V^2_{\max}}}
\right)^2-\frac{1}{V_{\max}^2}}.$$
The function $\gamma:\{(t,q)\in\setR^+\times\setR^+\,|\, t>\tilde t_0(q)\}\mapsto\setC$ is implicitly defined as the only root of the function 
$${\cal
  F}(\gamma,q,t)=z\left(\frac{1}{{\mat V^+_{Ps}}^2(q)}+\gamma^2\right)^{1/2}+h\left(\frac{1}{{\mat V_{Pf}^+}^2(q)}+\gamma^2\right)^{1/2}+i\gamma x-t$$
whose real part is positive.\\
The function $\upsilon:E_1\cup E_2\mapsto\setC$ is implicitly defined as the only root of the function 
$${\cal
  F}(\upsilon,q,t)=z\left(\frac{1}{{\mat V^+_{Ps}}^2(q)}+\upsilon^2\right)^{1/2}+h\left(\frac{1}{{\mat V_{Pf}^+}^2(q)}+\upsilon^2\right)^{1/2}+i\upsilon x-t$$
such that $\Im m\left[\partial_t \upsilon(t,q)
\right]<0$, with
$$E_1=\left\{
(t,q)\in\setR^+\times\setR^+\,|\, t_{h_1}<t<t_{0} \hbox{ and } 0<q<q_0(t)
\right\} $$ and  $$E_2=\left\{
(t,q)\in\setR^+\times\setR^+\,|\, t_{0}<t<t_{h_1} \hbox{ and } q_0(t)<q<q_1(t)\right\}.$$
  \item  $\gr u^+_{PfS}$ is the solid displacement of the
    reflected $PfS$ wave and satisfies:
$$
\left\{
\begin{array}{ll}
 \dsp u^+_{PfS,x}(x,z,t)=\frac{
   F^+_{Pf}}{\pi^2}\int_{0}^{q_{1}(t)}
\Re
 e\left[\ic\upsilon(t,q)\kapsip(\upsilon(t,q))\rfpsi(\upsilon(t,q))\frac{\partial \upsilon}{\partial t}(t,q)\right]dq,
\\[18pt]
\dsp u^+_{PfS,z}(x,z,t)= 
\frac{F^+_{Pf}}{\pi^2}\int_{0}^{q_{1}(t)}
\Re
e\left[(\upsilon^2(t,q)+q^2)\rfpsi(\upsilon(t,q))\frac{\partial\upsilon}{\partial    t}(t,q)\right]dq,
\end{array}\right.$$
$\dsp\hbox{if }t_{h_1} <t\leq t_{0} \hbox{ and } \left|\Im
    m\left[\gamma(t_{0},0)\right]\right|
  <\frac{1}{V_{\max}},$
$$\left\{
\begin{array}{lcl}
  \dsp u^+_{PfS,x}(x,z,t)&=&\dsp 
\begin{array}[t]{ll}
&\dsp\frac{\mat
    F^+_{Pf}}{\pi^2}
\int_{0}^{q_0(t)}
\Re
e\left[\ic\gamma(t,q)\kapsip(\gamma(t,q))\rfpsi(\gamma(t,q))\frac{\partial
    \gamma}{\partial t}(t,q)\right]dq
\\[18pt]
+&\dsp\frac{\mat
    F^+_{Pf}}{\pi^2}
\int_{q_0(t)}^{q_{1}(t)}
\Re
 e\left[\ic\upsilon(t,q)\kapsip(\upsilon(t,q))\rfpsi(\upsilon(t,q))\frac{\partial \upsilon}{\partial t}(t,q)\right]dq,
\end{array}
\\[60pt]
\dsp u^+_{PfS,z}(x,z,t)&=&\dsp 
\begin{array}[t]{ll}
&\dsp \frac{\mat
  F^+_{Pf}}{\pi^2}\int_{0}^{q_0(t)}
\Re
e\left[(\gamma^2(t,q)+q^2)\rfpsi(\gamma(t,q))\frac{\partial
    \gamma
}{\partial t}(t,q)\right]dq
\\[18pt]
+&\dsp\frac{\mat
  F^+_{Pf}}{\pi^2}
\int_{q_0(t)}^{q_{1}(t)}
\Re
e\left[(\upsilon^2(t,q)+q^2)\rfpsi(\upsilon(t,q))\frac{\partial\upsilon}{\partial    t}(t,q)\right]dq,
\end{array}
\end{array}\right.$$
$\dsp\hbox{if }t_{0} <t\leq t_{h_2}  \hbox{ and }  \left|\Im
    m\left[\gamma(t_{0},0)\right]\right|
  <\frac{1}{V_{\max}}$,
$$
\left\{
\begin{array}{ll}
  \dsp u^+_{PfS,x}(x,z,t)= \frac{F^+_{Pf}}{\pi^2}\int_{0}^{q_0(t)}
\Re
e\left[\ic\gamma(t,q)\kapsip(\gamma(t,q))\rfpsi(\gamma(t,q))\frac{\partial
    \gamma}{\partial t}(t,q)\right]dq,
\\[18pt]
\dsp u^+_{PfS,z}(x,z,t)=\frac{
 F^+_{Pf}}{\pi^2}\int_{0}^{q_0(t)}
\Re
e\left[(\gamma^2(t,q)+q^2)\rfpsi(\gamma(t,q))\frac{\partial
    \gamma
}{\partial t}(t,q)\right]dq,
\end{array}\right.$$
$\dsp\hbox{if } t_{h_2}<t \hbox{ and } \left|\Im
  m\left[\gamma(t_{0},0)\right]\right|
<\frac{1}{V_{\max}}$ or $\dsp\hbox{if }t_{0} <t
\hbox{ and } \left|\Im
  m\left[\gamma(t_{0},0)\right]\right|
\geq\frac{1}{V_{\max}}$\\
and $\gr
u^+_{PfS}(x,z,t)=0$ else.\\\\
$t_0$ denotes here the arrival time of the
reflected $PfS$ volume wave at point
$(x,0,z)$ (its calculation is similar to the
calculation of the arrival time of the transmitted
wave, see the appendix of~\cite{RAP2}),
\begin{eqnarray}
  t_{h_1}&=&h\sqrt{\frac{1}{{V^+_{Pf}}^2}-\frac{1}{V^2_{\max}}}+z\sqrt{\frac{1}{{V^+_{S}}^2}-\frac{1}{V^2_{\max}}}+\frac{|x|}{V_{\max}}
\end{eqnarray}
denotes the arrival time of the reflected $PfS$
head wave at point $(x,0,z)$,
\begin{eqnarray}
  t_{h_2}&=&\frac{h^2+z^2+hz\left(\frac{c_2}{c_1}+
\frac{c_1}{c_2}
\right)+x^2}
{\frac{h}{c_1}+\frac{z}{c_2}}
\end{eqnarray}
 denotes the time after which there is no longer head
wave at point $(x,0,z)$, where
$$c_1=\sqrt{\frac{1}{{V_{Pf}^+}^2}-\frac{1}{{V^2_{\max}}}} \hbox{ and } c_2=\sqrt{\frac{1}{{V^+_{S}}^2}-\frac{1}{{V^2_{\max}}}}.$$
The function $q_0:[t_0\,;\,+\infty]\mapsto\setR^+$ is
the reciprocal function of $\tilde t_0
:\setR^+\mapsto:[t_0,+\infty]$, where $\widetilde
t_0(q)$ is the arrival time  at point $(x,0,z)$ of the fictitious 
transmitted $PfS$ volume wave, propagating at a velocity $\mat
V_{Pf}^+(q)$ from the source to the interface and at velocity $\mat
V^+_{S}(q)$ from the interface to point $(x,0,z)$ (we refer again to~\cite{RAP2} for details on its calculation).\\[5pt]
The function $q_1:[t_1\,;\,t_0]\mapsto\setR^+$ is defined by 
$$q_{1}(t)=\sqrt{\frac{1}{x^2}\left(t-z\sqrt{\frac{1}{{V^+_{S}}^2}-\frac{1}{V^2_{\max}}}
-h\sqrt{\frac{1}{{V_{Pf}^+}^2}-\frac{1}{V^2_{\max}}}
\right)^2-\frac{1}{V_{\max}^2}}.$$
The function $\gamma:\{(t,q)\in\setR^+\times\setR^+\,|\, t>\tilde t_0(q)\}\mapsto\setC$ is implicitly defined as the only root of the function 
$${\cal
  F}(\gamma,q,t)=z\left(\frac{1}{{\mat V^+_{S}}^2(q)}+\gamma^2\right)^{1/2}+h\left(\frac{1}{{\mat V_{Pf}^+}^2(q)}+\gamma^2\right)^{1/2}+i\gamma x-t$$
whose real part is positive.\\
The function $\upsilon:E_1\cup E_2\mapsto\setC$ is implicitly defined as the only root of the function 
$${\cal
  F}(\upsilon,q,t)=z\left(\frac{1}{{\mat V^+_{S}}^2(q)}+\upsilon^2\right)^{1/2}+h\left(\frac{1}{{\mat V_{Pf}^+}^2(q)}+\upsilon^2\right)^{1/2}+i\upsilon x-t$$
such that $\Im m\left[\partial_t \upsilon(t,q)
\right]<0$, with
$$E_1=\left\{
(t,q)\in\setR^+\times\setR^+\,|\, t_{h_1}<t<t_{0} \hbox{ and } 0<q<q_0(t)
\right\} $$ and  $$E_2=\left\{
(t,q)\in\setR^+\times\setR^+\,|\, t_{0}<t<t_{h_1} \hbox{ and } q_0(t)<q<q_1(t)\right\}.$$
  \item  $\gr u^+_{PsPf}$ is the solid displacement of the
    reflected $PsPf$ wave and satisfies:
$$
\left\{
\begin{array}{ll}
 \dsp u^+_{PsPf,x}(x,z,t)=-\frac{\mat
   P_{11}^+F^+_{Ps}}{\pi^2}\int_{0}^{q_{1}(t)}
\Re
 e\left[\ic\upsilon(t,q)\rsf(\upsilon(t,q))\frac{\partial \upsilon}{\partial t}(t,q)\right]dq,
\\[18pt]
\dsp u^+_{PsPf,z}(x,z,t)= 
-\frac{\mat P_{11}^+F^+_{Ps}}{\pi^2}\int_{0}^{q_{1}(t)}
\Re
e\left[\kafp(\upsilon(t,q))\rsf(\upsilon(t,q))\frac{\partial\upsilon}{\partial    t}(t,q)\right]dq,
\end{array}\right.$$
$\dsp\hbox{if }t_{h_1} <t\leq t_{0} \hbox{ and } \left|\Im
    m\left[\gamma(t_{0},0)\right]\right|
  <\frac{1}{V_{\max}},$
$$\left\{
\begin{array}{lcl}
  \dsp u^+_{PsPf,x}(x,z,t)&=&\dsp 
\begin{array}[t]{ll}
-&\dsp\frac{\mat
    P_{11}^+F^+_{Ps}}{\pi^2}
\int_{0}^{q_0(t)}
\Re
e\left[\ic\gamma(t,q)\rsf(\gamma(t,q))\frac{\partial
    \gamma}{\partial t}(t,q)\right]dq
\\[18pt]
-&\dsp\frac{\mat
    P_{11}^+F^+_{Ps}}{\pi^2}
\int_{q_0(t)}^{q_{1}(t)}
\Re
 e\left[\ic\upsilon(t,q)\rsf(\upsilon(t,q))\frac{\partial \upsilon}{\partial t}(t,q)\right]dq,
\end{array}
\\[60pt]
\dsp u^+_{PsPf,z}(x,z,t)&=&\dsp 
\begin{array}[t]{ll}
-&\dsp \frac{\mat
  P_{11}^+F^+_{Ps}}{\pi^2}\int_{0}^{q_0(t)}
\Re
e\left[\kafp(\gamma(t,q))\rsf(\gamma(t,q))\frac{\partial
    \gamma
}{\partial t}(t,q)\right]dq
\\[18pt]
-&\dsp\frac{\mat
  P_{11}^+F^+_{Ps}}{\pi^2}
\int_{q_0(t)}^{q_{1}(t)}
\Re
e\left[\kafp(\upsilon(t,q))\rsf(\upsilon(t,q))\frac{\partial\upsilon}{\partial    t}(t,q)\right]dq,
\end{array}
\end{array}\right.$$
$\dsp\hbox{if }t_{0} <t\leq t_{h_2}  \hbox{ and }  \left|\Im
    m\left[\gamma(t_{0},0)\right]\right|
  <\frac{1}{V_{\max}}$,
$$
\left\{
\begin{array}{ll}
  \dsp u^+_{PsPf,x}(x,z,t)=- \frac{\mat P_{11}^+F^+_{Ps}}{\pi^2}\int_{0}^{q_0(t)}
\Re
e\left[\ic\gamma(t,q)\rsf(\gamma(t,q))\frac{\partial
    \gamma}{\partial t}(t,q)\right]dq,
\\[18pt]
\dsp u^+_{PsPf,z}(x,z,t)=-\frac{\mat
  P_{11}^+F^+_{Ps}}{\pi^2}\int_{0}^{q_0(t)}
\Re
e\left[\kafp(\gamma(t,q))\rsf(\gamma(t,q))\frac{\partial
    \gamma
}{\partial t}(t,q)\right]dq,
\end{array}\right.$$
$\dsp\hbox{if } t_{h_2}<t \hbox{ and } \left|\Im
  m\left[\gamma(t_{0},0)\right]\right|
<\frac{1}{V_{\max}}$ or $\dsp\hbox{if }t_{0} <t
\hbox{ and } \left|\Im
  m\left[\gamma(t_{0},0)\right]\right|
\geq\frac{1}{V_{\max}}$\\
and $\gr
u^+_{PsPf}(x,z,t)=0$ else.\\\\
$t_0$ denotes here the arrival time of the
reflected $PsPf$ volume wave at point
$(x,0,z)$,
\begin{eqnarray}
  t_{h_1}&=&h\sqrt{\frac{1}{{V^+_{Ps}}^2}-\frac{1}{V^2_{\max}}}+z\sqrt{\frac{1}{{V^+_{Pf}}^2}-\frac{1}{V^2_{\max}}}+\frac{|x|}{V_{\max}}
\end{eqnarray}
denotes the arrival time of the reflected $PsPf$
head wave at point $(x,0,z)$,
\begin{eqnarray}
  t_{h_2}&=&\frac{h^2+z^2+hz\left(\frac{c_2}{c_1}+
\frac{c_1}{c_2}
\right)+x^2}
{\frac{h}{c_1}+\frac{z}{c_2}}
\end{eqnarray}
 denotes the time after which there is no longer head
wave at point $(x,0,z)$, where
$$c_1=\sqrt{\frac{1}{{V_{Ps}^+}^2}-\frac{1}{{V^2_{\max}}}} \hbox{ and } c_2=\sqrt{\frac{1}{{V^+_{Pf}}^2}-\frac{1}{{V^2_{\max}}}}.$$
The function $q_0:[t_0\,;\,+\infty]\mapsto\setR^+$ is
the reciprocal function of $\tilde t_0
:\setR^+\mapsto:[t_0,+\infty]$, where $\widetilde
t_0(q)$ is the arrival time  at point $(x,0,z)$ of the fictitious 
transmitted $PsPf$ volume wave, propagating at a velocity $\mat
V_{Ps}^+(q)$ from the source to the interface and at velocity $\mat
V^+_{Pf}(q)$ from the interface to point $(x,0,z)$.\\[5pt]
The function $q_1:[t_1\,;\,t_0]\mapsto\setR^+$ is defined by 
$$q_{1}(t)=\sqrt{\frac{1}{x^2}\left(t-z\sqrt{\frac{1}{{V^+_{Pf}}^2}-\frac{1}{V^2_{\max}}}
-h\sqrt{\frac{1}{{V_{Ps}^+}^2}-\frac{1}{V^2_{\max}}}
\right)^2-\frac{1}{V_{\max}^2}}.$$
The function $\gamma:\{(t,q)\in\setR^+\times\setR^+\,|\, t>\tilde t_0(q)\}\mapsto\setC$ is implicitly defined as the only root of the function 
$${\cal
  F}(\gamma,q,t)=z\left(\frac{1}{{\mat V^+_{Pf}}^2(q)}+\gamma^2\right)^{1/2}+h\left(\frac{1}{{\mat V_{Ps}^+}^2(q)}+\gamma^2\right)^{1/2}+i\gamma x-t$$
whose real part is positive.\\
The function $\upsilon:E_1\cup E_2\mapsto\setC$ is implicitly defined as the only root of the function 
$${\cal
  F}(\upsilon,q,t)=z\left(\frac{1}{{\mat V^+_{Pf}}^2(q)}+\upsilon^2\right)^{1/2}+h\left(\frac{1}{{\mat V_{Ps}^+}^2(q)}+\upsilon^2\right)^{1/2}+i\upsilon x-t$$
such that $\Im m\left[\partial_t \upsilon(t,q)
\right]<0$, with
$$E_1=\left\{
(t,q)\in\setR^+\times\setR^+\,|\, t_{h_1}<t<t_{0} \hbox{ and } 0<q<q_0(t)
\right\} $$ and  $$E_2=\left\{
(t,q)\in\setR^+\times\setR^+\,|\, t_{0}<t<t_{h_1} \hbox{ and } q_0(t)<q<q_1(t)\right\}.$$
  \item  $\gr u^+_{PsPs}$ is the solid displacement of the
    reflected $PsPs$ wave and satisfies:
$$\left\{
\begin{array}{lcl}
 \dsp  u_{PsPs,x}^+(x,z,t)&=&\dsp\mat P^+_{12}F^+_{Ps}\int_{0}^{q_1(t)} \frac{\Im
   m\Big[\ic\upsilon(t,q)\kasp(\upsilon(t,q))\rss(\upsilon(t,q))\Big]}{\pi^2r \sqrt{q^2+q_0^2(t)}}dq,
\\[15pt]
\dsp u_{PsPs,z}^+(x,z,t)&=&\dsp \mat P^+_{12}F^+_{Ps}
\int_{0}^{q_1(t)}\frac{\Im
  m\Big[\kasp^2(\upsilon(t,q))\rss(\upsilon(t,q))\Big]}{\pi^2r \sqrt{q^2+q_0^2(t)}}dq,
\end{array}\right.$$
$\dsp\hbox{if }t_{\hbox{h}_1} <t\leq t_{0} \hbox{ and } \frac{x}{r}>\frac{V_{Ps}^+}{V_{\max}}$,
$$\left\{
\begin{array}{lcl}
  \dsp u_{PsPs,x}^+(x,z,t)&=&\dsp 
 \begin{array}[t]{ll}
&\dsp \mat P^+_{12}F^+_{Ps} \int^{q_1(t)}_{q_0(t)}\frac{\Im
   m\Big[(\ic\upsilon(t,q)\kasp(\upsilon(t,q))\rss(\upsilon(t,q))\Big]}{\pi^2r \sqrt{q^2-q_0^2(t)}}dq
\\[18pt]
 -&\dsp\mat P^+_{12}F^+_{Ps}\int_{0}^{q_0(t)}\frac{\Re
   e\Big[\ic\gamma(t,q)\kasp(\gamma(t,q))\rss(\gamma(t,q))\Big]}{\pi^2r \sqrt{q_0^2(t)-q^2}}dq,
\end{array}
\\[60pt]
\dsp u_{PsPs,y}^+(x,z,t)&=&\dsp 
  \begin{array}[t]{ll}
&\dsp\mat P^+_{12}F^+_{Ps}\int^{q_1(t)}_{q_0(t)}\frac{\Im
  m\Big[\kasp^2(\upsilon(t,q))\rss(\upsilon(t,q))\Big]}{\pi^2r \sqrt{q^2-q_0^2(t)}}dq
\\[15pt]
-&\dsp \mat P^+_{12}F^+_{Ps}\int_{0}^{q_0(t)} \frac{\Re
  e\Big[\kasp^2(\gamma(t,q))\rss(\gamma(t,q))\Big]}{\pi^2r \sqrt{q_0^2(t)-q^2}}dq,
\end{array}
\end{array}\right.$$
$\dsp\hbox{if }t_{0} <t\leq t_{\hbox{h}_2}  \hbox{ and } \frac{x}{r}>\frac{V_{Ps}^+}{V_{\max}}$,
$$\left\{
\begin{array}{lcl}
  \dsp u_{PsPs,x}^+(x,z,t)&=&-\mat P^+_{12}F^+_{Ps}\dsp  \int^{q_0(t)}_{0}\frac{\Re e\Big[\ic\gamma(t,q)\kasp(\gamma(t,q))\rss(\gamma(t,q))\Big]}{\pi^2r \sqrt{q_0^2(t)-q^2}}dq,
\\[15pt]
\dsp u_{PsPs,y}^+(x,z,t)&=&-\mat P^+_{12}F^+_{Ps}\dsp \int^{q_0(t)}_{0} \frac{\Re e\Big[\kasp^2(\gamma(t,q))\rss(\gamma(t,q))\Big]}{\pi^2r \sqrt{q_0^2(t)-q^2}}dq,
\end{array}\right.$$
$\dsp\hbox{if } t_{\hbox{h}_2}<t
\hbox{ and } \frac{x}{r}>\frac{V_{Ps}^+}{V_{\max}}$ or
 $\dsp\hbox{if }t_{0} <t
\hbox{ and } \frac{x}{r}\leq \frac{V_{Ps}^+}{V_{\max}}$
and 
$$
 \gr
u_{PsPs}(x,z,t)=0 \hbox{ else }.$$
We set here $r=(x^2+(z+h)^2)^{1/2}$  and
$t_0=r/V_{Ps}^+$ denotes the arrival time of the reflected $PsPs$
volume wave at point $(x,0,z)$, 
\begin{eqnarray}
  t_{h_1}&=&(z+h)\sqrt{\frac{1}{{V_{Ps}^+}^2}-\frac{1}{V^2_{\max}}}+\frac{|x|}{V_{\max}}
\end{eqnarray}
denotes the arrival time of the reflected $PsPs$
head-wave at point $(x,0,z)$ and 
\begin{eqnarray}
  t_{h_2}&=&\frac{r}{z+h}\sqrt{\frac{1}{{V_{Ps}^+}^2}-\frac{1}{{V^2_{\max}}}}
\end{eqnarray}
denotes the time after which there is no longer head
wave at point $(x,0,z)$, (contrary to the 2D case, this time does not
coincide with the arrival time of the volume wave).
We also define the functions $\gamma$, $\upsilon$,
 $q_0$ and $q_1$ by
$$\gamma : \{t\in\setR \,|\, t>t_0\}\times\setR\mapsto
\setC:=
\gamma(t,q_y)=\ic \frac
{xt}{r^2}+\frac{z+h}{r}\sqrt{\frac{t^2}{r^2}-\frac{1}{{\mat
      V_{Ps}^{+2}}(q_y)}}$$
$$\upsilon : \{t\in\setR \,|\, t_{h_1}<t<t_{h_2}\}\times\setR\mapsto
\setC:=
\upsilon(t,q_y)=-\ic\left(\frac{z+h}{r}-\sqrt{\frac{1}{\mat
      V_{Ps}^{+2}(q_y)}-\frac{t^2}{r^2}}+\frac
  {x}{r^2}t\right),$$
$$q_0:\setR\to\setR:=q_0(t)=\sqrt{\left|\frac{t^2}{r^2}-\frac{1}{
      {V_{Ps}^+}^2}\right|}$$
and 
$$q_1:\setR\to\setR:=q_1(t)=\sqrt{\frac{1}{x^2}\left(t-(z+h)\sqrt{\frac{1}{{V_{Ps}^+}^2}-\frac{1}{V^2_{\max}}}\right)^2-\frac{1}{V_{\max}^2}}.$$
  \item  $\gr u^+_{PsS}$ is the solid displacement of the
    reflected $PsS$ wave and satisfies:
$$
\left\{
\begin{array}{ll}
 \dsp u^+_{PsS,x}(x,z,t)=\frac{
   F^+_{Ps}}{\pi^2}\int_{0}^{q_{1}(t)}
\Re
 e\left[\ic\upsilon(t,q)\kapsip(\upsilon(t,q))\rspsi(\upsilon(t,q))\frac{\partial \upsilon}{\partial t}(t,q)\right]dq,
\\[18pt]
\dsp u^+_{PsS,z}(x,z,t)= 
\frac{F^+_{Ps}}{\pi^2}\int_{0}^{q_{1}(t)}
\Re
e\left[(\upsilon^2(t,q)+q^2)\rspsi(\upsilon(t,q))\frac{\partial\upsilon}{\partial    t}(t,q)\right]dq,
\end{array}\right.$$
$\dsp\hbox{if }t_{h_1} <t\leq t_{0} \hbox{ and } \left|\Im
    m\left[\gamma(t_{0},0)\right]\right|
  <\frac{1}{V_{\max}},$
$$\left\{
\begin{array}{lcl}
  \dsp u^+_{PsS,x}(x,z,t)&=&\dsp 
\begin{array}[t]{ll}
&\dsp\frac{\mat
    F^+_{Ps}}{\pi^2}
\int_{0}^{q_0(t)}
\Re
e\left[\ic\gamma(t,q)\kapsip(\gamma(t,q))\rspsi(\gamma(t,q))\frac{\partial
    \gamma}{\partial t}(t,q)\right]dq
\\[18pt]
+&\dsp\frac{\mat
    F^+_{Ps}}{\pi^2}
\int_{q_0(t)}^{q_{1}(t)}
\Re
 e\left[\ic\upsilon(t,q)\kapsip(\upsilon(t,q))\rspsi(\upsilon(t,q))\frac{\partial \upsilon}{\partial t}(t,q)\right]dq,
\end{array}
\\[60pt]
\dsp u^+_{PsS,z}(x,z,t)&=&\dsp 
\begin{array}[t]{ll}
&\dsp \frac{\mat
  F^+_{Ps}}{\pi^2}\int_{0}^{q_0(t)}
\Re
e\left[(\gamma^2(t,q)+q^2)\rspsi(\gamma(t,q))\frac{\partial
    \gamma
}{\partial t}(t,q)\right]dq
\\[18pt]
+&\dsp\frac{\mat
  F^+_{Ps}}{\pi^2}
\int_{q_0(t)}^{q_{1}(t)}
\Re
e\left[(\upsilon^2(t,q)+q^2)\rspsi(\upsilon(t,q))\frac{\partial\upsilon}{\partial    t}(t,q)\right]dq,
\end{array}
\end{array}\right.$$
$\dsp\hbox{if }t_{0} <t\leq t_{h_2}  \hbox{ and }  \left|\Im
    m\left[\gamma(t_{0},0)\right]\right|
  <\frac{1}{V_{\max}}$,
$$
\left\{
\begin{array}{ll}
  \dsp u^+_{PsS,x}(x,z,t)= \frac{F^+_{Ps}}{\pi^2}\int_{0}^{q_0(t)}
\Re
e\left[\ic\gamma(t,q)\kapsip(\gamma(t,q))\rspsi(\gamma(t,q))\frac{\partial
    \gamma}{\partial t}(t,q)\right]dq,
\\[18pt]
\dsp u^+_{PsS,z}(x,z,t)=\frac{
 F^+_{Ps}}{\pi^2}\int_{0}^{q_0(t)}
\Re
e\left[(\gamma^2(t,q)+q^2)\rspsi(\gamma(t,q))\frac{\partial
    \gamma
}{\partial t}(t,q)\right]dq,
\end{array}\right.$$
$\dsp\hbox{if } t_{h_2}<t \hbox{ and } \left|\Im
  m\left[\gamma(t_{0},0)\right]\right|
<\frac{1}{V_{\max}}$ or $\dsp\hbox{if }t_{0} <t
\hbox{ and } \left|\Im
  m\left[\gamma(t_{0},0)\right]\right|
\geq\frac{1}{V_{\max}}$\\
and $\gr
u^+_{PsS}(x,z,t)=0$ else.\\\\
$t_0$ denotes here the arrival time of the
reflected $PsS$ volume wave at point
$(x,0,z)$ (its calculation is similar to the
calculation of the arrival time of the transmitted
wave, see the appendix of~\cite{RAP2}),
\begin{eqnarray}
  t_{h_1}&=&h\sqrt{\frac{1}{{V^+_{Ps}}^2}-\frac{1}{V^2_{\max}}}+z\sqrt{\frac{1}{{V^+_{S}}^2}-\frac{1}{V^2_{\max}}}+\frac{|x|}{V_{\max}}
\end{eqnarray}
denotes the arrival time of the reflected $PsS$
head wave at point $(x,0,z)$,
\begin{eqnarray}
  t_{h_2}&=&\frac{h^2+z^2+hz\left(\frac{c_2}{c_1}+
\frac{c_1}{c_2}
\right)+x^2}
{\frac{h}{c_1}+\frac{z}{c_2}}
\end{eqnarray}
 denotes the time after which there is no longer head
wave at point $(x,0,z)$, where
$$c_1=\sqrt{\frac{1}{{V_{Ps}^+}^2}-\frac{1}{{V^2_{\max}}}} \hbox{ and } c_2=\sqrt{\frac{1}{{V^+_{S}}^2}-\frac{1}{{V^2_{\max}}}}.$$
The function $q_0:[t_0\,;\,+\infty]\mapsto\setR^+$ is
the reciprocal function of $\tilde t_0
:\setR^+\mapsto:[t_0,+\infty]$, where $\widetilde
t_0(q)$ is the arrival time  at point $(x,0,z)$ of the fictitious 
transmitted $PsS$ volume wave, propagating at a velocity $\mat
V_{Ps}^+(q)$ from the source to the interface and at velocity $\mat
V^+_{S}(q)$ from the interface to point $(x,0,z)$ (we refer again to~\cite{RAP2} for details on its calculation).\\[5pt]
The function $q_1:[t_1\,;\,t_0]\mapsto\setR^+$ is defined by 
$$q_{1}(t)=\sqrt{\frac{1}{x^2}\left(t-z\sqrt{\frac{1}{{V^+_{S}}^2}-\frac{1}{V^2_{\max}}}
-h\sqrt{\frac{1}{{V_{Ps}^+}^2}-\frac{1}{V^2_{\max}}}
\right)^2-\frac{1}{V_{\max}^2}}.$$
The function $\gamma:\{(t,q)\in\setR^+\times\setR^+\,|\, t>\tilde t_0(q)\}\mapsto\setC$ is implicitly defined as the only root of the function 
$${\cal
  F}(\gamma,q,t)=z\left(\frac{1}{{\mat V^+_{S}}^2(q)}+\gamma^2\right)^{1/2}+h\left(\frac{1}{{\mat V_{Ps}^+}^2(q)}+\gamma^2\right)^{1/2}+i\gamma x-t$$
whose real part is positive.\\
The function $\upsilon:E_1\cup E_2\mapsto\setC$ is implicitly defined as the only root of the function 
$${\cal
  F}(\upsilon,q,t)=z\left(\frac{1}{{\mat V^+_{S}}^2(q)}+\upsilon^2\right)^{1/2}+h\left(\frac{1}{{\mat V_{Ps}^+}^2(q)}+\upsilon^2\right)^{1/2}+i\upsilon x-t$$
such that $\Im m\left[\partial_t \upsilon(t,q)
\right]<0$, with
$$E_1=\left\{
(t,q)\in\setR^+\times\setR^+\,|\, t_{h_1}<t<t_{0} \hbox{ and } 0<q<q_0(t)
\right\} $$ and  $$E_2=\left\{
(t,q)\in\setR^+\times\setR^+\,|\, t_{0}<t<t_{h_1} \hbox{ and } q_0(t)<q<q_1(t)\right\}.$$
  \item  $\gr u^-_{PfPf}$ is the solid displacement of the
    transmitted $PfPf$ wave (the $Pf$ transmitted wave generated by the $Pf$ incident wave) and satisfies:
$$
\left\{
\begin{array}{ll}
 \dsp u^-_{PfPf,x}(x,z,t)=-\frac{\mat
   P_{11}^-F^+_{Pf}}{\pi^2}\int_{0}^{q_{1}(t)}
\Re
 e\left[\ic\upsilon(t,q)\tff(\upsilon(t,q))\frac{\partial \upsilon}{\partial t}(t,q)\right]dq,
\\[18pt]
\dsp u^-_{PfPf,z}(x,z,t)= 
\frac{\mat P_{11}^-F^+_{Pf}}{\pi^2}\int_{0}^{q_{1}(t)}
\Re
e\left[\kaf(\upsilon(t,q))\tff(\upsilon(t,q))\frac{\partial\upsilon}{\partial    t}(t,q)\right]dq,
\end{array}\right.$$
$\dsp\hbox{if }t_{h_1} <t\leq t_{0} \hbox{ and } \left|\Im
    m\left[\gamma(t_{0},0)\right]\right|
  <\frac{1}{V_{\max}},$
$$\left\{
\begin{array}{lcl}
  \dsp u^-_{PfPf,x}(x,z,t)&=&\dsp 
\begin{array}[t]{ll}
-&\dsp\frac{\mat
    P_{11}^-F^+_{Pf}}{\pi^2}
\int_{0}^{q_0(t)}
\Re
e\left[\ic\gamma(t,q)\tff(\gamma(t,q))\frac{\partial
    \gamma}{\partial t}(t,q)\right]dq
\\[18pt]
-&\dsp\frac{\mat
    P_{11}^-F^+_{Pf}}{\pi^2}
\int_{q_0(t)}^{q_{1}(t)}
\Re
 e\left[\ic\upsilon(t,q)\tff(\upsilon(t,q))\frac{\partial \upsilon}{\partial t}(t,q)\right]dq,
\end{array}
\\[60pt]
\dsp u^-_{PfPf,z}(x,z,t)&=&\dsp 
\begin{array}[t]{ll}
&\dsp \frac{\mat
  P_{11}^-F^+_{Pf}}{\pi^2}\int_{0}^{q_0(t)}
\Re
e\left[\kaf(\gamma(t,q))\tff(\gamma(t,q))\frac{\partial
    \gamma
}{\partial t}(t,q)\right]dq
\\[18pt]
+&\dsp\frac{\mat
  P_{11}^-F^+_{Pf}}{\pi^2}
\int_{q_0(t)}^{q_{1}(t)}
\Re
e\left[\kaf(\upsilon(t,q))\tff(\upsilon(t,q))\frac{\partial\upsilon}{\partial    t}(t,q)\right]dq,
\end{array}
\end{array}\right.$$
$\dsp\hbox{if }t_{0} <t\leq t_{h_2}  \hbox{ and }  \left|\Im
    m\left[\gamma(t_{0},0)\right]\right|
  <\frac{1}{V_{\max}}$,
$$
\left\{
\begin{array}{ll}
  \dsp u^-_{PfPf,x}(x,z,t)=- \frac{\mat P_{11}^-F^+_{Pf}}{\pi^2}\int_{0}^{q_0(t)}
\Re
e\left[\ic\gamma(t,q)\tff(\gamma(t,q))\frac{\partial
    \gamma}{\partial t}(t,q)\right]dq,
\\[18pt]
\dsp u^-_{PfPf,z}(x,z,t)=\frac{\mat
  P_{11}^-F^+_{Pf}}{\pi^2}\int_{0}^{q_0(t)}
\Re
e\left[\kaf(\gamma(t,q))\tff(\gamma(t,q))\frac{\partial
    \gamma
}{\partial t}(t,q)\right]dq,
\end{array}\right.$$
$\dsp\hbox{if } t_{h_2}<t
\hbox{ and } \left|\Im
    m\left[\gamma(t_{0},0)\right]\right|
  <\frac{1}{V_{\max}}$ or
 $\dsp\hbox{if }t_{0} <t
\hbox{ and } \left|\Im
    m\left[\gamma(t_{0},0)\right]\right|
  \geq\frac{1}{V_{\max}}$\\
and  $\gr
 u^-_{PfPf}(x,z,t)=0$ else.\\\\
 $t_0$ denotes here the arrival time of the transmitted $PfPf$
 volume wave at point $(x,0,z)$,
\begin{eqnarray}
  t_{h_1}&=&h\sqrt{\frac{1}{{V^+_{Pf}}^2}-\frac{1}{V^2_{\max}}}-z\sqrt{\frac{1}{{V^-_{Pf}}^2}-\frac{1}{V^2_{\max}}}+\frac{|x|}{V_{\max}}
\end{eqnarray}
denotes the arrival time of the transmitted $PfPf$
head wave at point $(x,0,z)$,
\begin{eqnarray}
  t_{h_2}&=&\frac{h^2+z^2-hz\left(\frac{c_2}{c_1}+
\frac{c_1}{c_2}
\right)+x^2}
{\frac{h}{c_1}-\frac{z}{c_2}}
\end{eqnarray}
 denotes the time after which there is no longer head
wave at point $(x,0,z)$, where
$$c_1=\sqrt{\frac{1}{{V_{Pf}^+}^2}-\frac{1}{{V^2_{\max}}}} \hbox{ and } c_2=\sqrt{\frac{1}{{V^-_{Pf}}^2}-\frac{1}{{V^2_{\max}}}}.$$
The function $q_0:[t_0\,;\,+\infty]\mapsto\setR^+$ is
the reciprocal function of $\tilde t_0
:\setR^+\mapsto:[t_0,+\infty]$, where $\widetilde
t_0(q)$ is the arrival time  at point $(x,0,z)$ of the fictitious 
transmitted $PfPf$ volume wave, propagating at a velocity $\mat
V_{Pf}^+(q)$ in the top layer and at velocity $\mat
V^-_{Pf}(q)$ in the bottom layer.\\[5pt]
The function $q_1:[t_1\,;\,t_0]\mapsto\setR^+$ is defined by 
$$q_{1}(t)=\sqrt{\frac{1}{x^2}\left(t+z\sqrt{\frac{1}{{V^-_{Pf}}^2}-\frac{1}{V^2_{\max}}}
-h\sqrt{\frac{1}{{V_{Pf}^+}^2}-\frac{1}{V^2_{\max}}}
\right)^2-\frac{1}{V_{\max}^2}}.$$
The function $\gamma:\{(t,q)\in\setR^+\times\setR^+\,|\, t>\tilde t_0(q)\}\mapsto\setC$ is implicitly defined as the only root of the function 
$${\cal
  F}(\gamma,q,t)=-z\left(\frac{1}{{\mat V^-_{Pf}}^2(q)}+\gamma^2\right)^{1/2}+h\left(\frac{1}{{\mat V_{Pf}^+}^2(q)}+\gamma^2\right)^{1/2}+i\gamma x-t$$
whose real part is positive.\\
The function $\upsilon:E_1\cup E_2\mapsto\setC$ is implicitly defined as the only root of the function 
$${\cal
  F}(\upsilon,q,t)=-z\left(\frac{1}{{\mat V^-_{Pf}}^2(q)}+\upsilon^2\right)^{1/2}+h\left(\frac{1}{{\mat V_{Pf}^+}^2(q)}+\upsilon^2\right)^{1/2}+i\upsilon x-t$$
such that $\Im m\left[\partial_t \upsilon(t,q)
\right]<0$, with
$$E_1=\left\{
(t,q)\in\setR^+\times\setR^+\,|\, t_{h_1}<t<t_{0} \hbox{ and } 0<q<q_0(t)
\right\} $$ and  $$E_2=\left\{
(t,q)\in\setR^+\times\setR^+\,|\, t_{0}<t<t_{h_1} \hbox{ and } q_0(t)<q<q_1(t)\right\}.$$
  \item  $\gr u^-_{PfPs}$ is the solid displacement of the
    transmitted $PfPs$ wave and satisfies:
$$
\left\{
\begin{array}{ll}
 \dsp u^-_{PfPs,x}(x,z,t)=-\frac{\mat
   P_{12}^-F^+_{Pf}}{\pi^2}\int_{0}^{q_{1}(t)}
\Re
 e\left[\ic\upsilon(t,q)\tfs(\upsilon(t,q))\frac{\partial \upsilon}{\partial t}(t,q)\right]dq,
\\[18pt]
\dsp u^-_{PfPs,z}(x,z,t)= 
\frac{\mat P_{12}^-F^+_{Pf}}{\pi^2}\int_{0}^{q_{1}(t)}
\Re
e\left[\kas(\upsilon(t,q))\tfs(\upsilon(t,q))\frac{\partial\upsilon}{\partial    t}(t,q)\right]dq,
\end{array}\right.$$
$\dsp\hbox{if }t_{h_1} <t\leq t_{0} \hbox{ and } \left|\Im
    m\left[\gamma(t_{0},0)\right]\right|
  <\frac{1}{V_{\max}},$
$$\left\{
\begin{array}{lcl}
  \dsp u^-_{PfPs,x}(x,z,t)&=&\dsp 
\begin{array}[t]{ll}
-&\dsp\frac{\mat
    P_{12}^-F^+_{Pf}}{\pi^2}
\int_{0}^{q_0(t)}
\Re
e\left[\ic\gamma(t,q)\tfs(\gamma(t,q))\frac{\partial
    \gamma}{\partial t}(t,q)\right]dq
\\[18pt]
-&\dsp\frac{\mat
    P_{12}^-F^+_{Pf}}{\pi^2}
\int_{q_0(t)}^{q_{1}(t)}
\Re
 e\left[\ic\upsilon(t,q)\tfs(\upsilon(t,q))\frac{\partial \upsilon}{\partial t}(t,q)\right]dq,
\end{array}
\\[60pt]
\dsp u^-_{PfPs,z}(x,z,t)&=&\dsp 
\begin{array}[t]{ll}
&\dsp \frac{\mat
  P_{12}^-F^+_{Pf}}{\pi^2}\int_{0}^{q_0(t)}
\Re
e\left[\kas(\gamma(t,q))\tfs(\gamma(t,q))\frac{\partial
    \gamma
}{\partial t}(t,q)\right]dq
\\[18pt]
+&\dsp\frac{\mat
  P_{12}^-F^+_{Pf}}{\pi^2}
\int_{q_0(t)}^{q_{1}(t)}
\Re
e\left[\kas(\upsilon(t,q))\tfs(\upsilon(t,q))\frac{\partial\upsilon}{\partial    t}(t,q)\right]dq,
\end{array}
\end{array}\right.$$
$\dsp\hbox{if }t_{0} <t\leq t_{h_2}  \hbox{ and }  \left|\Im
    m\left[\gamma(t_{0},0)\right]\right|
  <\frac{1}{V_{\max}}$,
$$
\left\{
\begin{array}{ll}
  \dsp u^-_{PfPs,x}(x,z,t)=- \frac{\mat P_{12}^-F^+_{Pf}}{\pi^2}\int_{0}^{q_0(t)}
\Re
e\left[\ic\gamma(t,q)\tfs(\gamma(t,q))\frac{\partial
    \gamma}{\partial t}(t,q)\right]dq,
\\[18pt]
\dsp u^-_{PfPs,z}(x,z,t)=\frac{\mat
  P_{12}^-F^+_{Pf}}{\pi^2}\int_{0}^{q_0(t)}
\Re
e\left[\kas(\gamma(t,q))\tfs(\gamma(t,q))\frac{\partial
    \gamma
}{\partial t}(t,q)\right]dq,
\end{array}\right.$$
$\dsp\hbox{if } t_{h_2}<t
\hbox{ and } \left|\Im
    m\left[\gamma(t_{0},0)\right]\right|
  <\frac{1}{V_{\max}}$ or
 $\dsp\hbox{if }t_{0} <t
\hbox{ and } \left|\Im
    m\left[\gamma(t_{0},0)\right]\right|
  \geq\frac{1}{V_{\max}}$\\
and  $\gr
 u^-_{PfPs}(x,z,t)=0$ else.\\\\
 $t_0$ denotes here the arrival time of the transmitted $PfPs$
 volume wave at point $(x,0,z)$,
\begin{eqnarray}
  t_{h_1}&=&h\sqrt{\frac{1}{{V^+_{Pf}}^2}-\frac{1}{V^2_{\max}}}-z\sqrt{\frac{1}{{V^-_{Ps}}^2}-\frac{1}{V^2_{\max}}}+\frac{|x|}{V_{\max}}
\end{eqnarray}
denotes the arrival time of the transmitted $PfPs$
head wave at point $(x,0,z)$,
\begin{eqnarray}
  t_{h_2}&=&\frac{h^2+z^2-hz\left(\frac{c_2}{c_1}+
\frac{c_1}{c_2}
\right)+x^2}
{\frac{h}{c_1}-\frac{z}{c_2}}
\end{eqnarray}
 denotes the time after which there is no longer head
wave at point $(x,0,z)$, where
$$c_1=\sqrt{\frac{1}{{V_{Pf}^+}^2}-\frac{1}{{V^2_{\max}}}} \hbox{ and } c_2=\sqrt{\frac{1}{{V^-_{Ps}}^2}-\frac{1}{{V^2_{\max}}}}.$$
The function $q_0:[t_0\,;\,+\infty]\mapsto\setR^+$ is
the reciprocal function of $\tilde t_0
:\setR^+\mapsto:[t_0,+\infty]$, where $\widetilde
t_0(q)$ is the arrival time  at point $(x,0,z)$ of the fictitious 
transmitted $PfPs$ volume wave, propagating at a velocity $\mat
V_{Pf}^+(q)$ in the top layer and at velocity $\mat
V^-_{Ps}(q)$ in the bottom layer.\\[5pt]
The function $q_1:[t_1\,;\,t_0]\mapsto\setR^+$ is defined by 
$$q_{1}(t)=\sqrt{\frac{1}{x^2}\left(t+z\sqrt{\frac{1}{{V^-_{Ps}}^2}-\frac{1}{V^2_{\max}}}
-h\sqrt{\frac{1}{{V_{Pf}^+}^2}-\frac{1}{V^2_{\max}}}
\right)^2-\frac{1}{V_{\max}^2}}.$$
The function $\gamma:\{(t,q)\in\setR^+\times\setR^+\,|\, t>\tilde t_0(q)\}\mapsto\setC$ is implicitly defined as the only root of the function 
$${\cal
  F}(\gamma,q,t)=-z\left(\frac{1}{{\mat V^-_{Ps}}^2(q)}+\gamma^2\right)^{1/2}+h\left(\frac{1}{{\mat V_{Pf}^+}^2(q)}+\gamma^2\right)^{1/2}+i\gamma x-t$$
whose real part is positive.\\
The function $\upsilon:E_1\cup E_2\mapsto\setC$ is implicitly defined as the only root of the function 
$${\cal
  F}(\upsilon,q,t)=-z\left(\frac{1}{{\mat V^-_{Ps}}^2(q)}+\upsilon^2\right)^{1/2}+h\left(\frac{1}{{\mat V_{Pf}^+}^2(q)}+\upsilon^2\right)^{1/2}+i\upsilon x-t$$
such that $\Im m\left[\partial_t \upsilon(t,q)
\right]<0$, with
$$E_1=\left\{
(t,q)\in\setR^+\times\setR^+\,|\, t_{h_1}<t<t_{0} \hbox{ and } 0<q<q_0(t)
\right\} $$ and  $$E_2=\left\{
(t,q)\in\setR^+\times\setR^+\,|\, t_{0}<t<t_{h_1} \hbox{ and } q_0(t)<q<q_1(t)\right\}.$$
  \item  $\gr u^-_{PfS}$ is the solid displacement of the
    transmitted $PfS$ wave and satisfies:
$$
\left\{
\begin{array}{ll}
 \dsp u^-_{PfS,x}(x,z,t)=-\frac{
   F^+_{Pf}}{\pi^2}\int_{0}^{q_{1}(t)}
\Re
 e\left[\ic\upsilon(t,q)\kapsi(\upsilon(t,q))\tfpsi(\upsilon(t,q))\frac{\partial \upsilon}{\partial t}(t,q)\right]dq,
\\[18pt]
\dsp u^-_{PfS,z}(x,z,t)= 
\frac{F^+_{Pf}}{\pi^2}\int_{0}^{q_{1}(t)}
\Re
e\left[(\upsilon^2(t,q)+q^2)\tfpsi(\upsilon(t,q))\frac{\partial\upsilon}{\partial    t}(t,q)\right]dq,
\end{array}\right.$$
$\dsp\hbox{if }t_{h_1} <t\leq t_{0} \hbox{ and } \left|\Im
    m\left[\gamma(t_{0},0)\right]\right|
  <\frac{1}{V_{\max}},$
$$\left\{
\begin{array}{lcl}
  \dsp u^-_{PfS,x}(x,z,t)&=&\dsp 
\begin{array}[t]{ll}
-&\dsp\frac{
  F^+_{Pf}}{\pi^2}
\int_{0}^{q_0(t)}
\Re
e\left[\ic\gamma(t,q)\kapsi(\gamma(t,q))\tfpsi(\gamma(t,q))\frac{\partial
    \gamma}{\partial t}(t,q)\right]dq
\\[18pt]
-&\dsp\frac{
F^+_{Pf}}{\pi^2}
\int_{q_0(t)}^{q_{1}(t)}
\Re
 e\left[\ic\upsilon(t,q)\kapsi(\upsilon(t,q))\tfpsi(\upsilon(t,q))\frac{\partial \upsilon}{\partial t}(t,q)\right]dq,
\end{array}
\\[60pt]
\dsp u^-_{PfS,z}(x,z,t)&=&\dsp 
\begin{array}[t]{ll}
&\dsp \frac{\mat
F^+_{Pf}}{\pi^2}\int_{0}^{q_0(t)}
\Re
e\left[(\gamma^2(t,q)+q^2)\tfpsi(\gamma(t,q))\frac{\partial
    \gamma
}{\partial t}(t,q)\right]dq
\\[18pt]
+&\dsp\frac{\mat
 F^+_{Pf}}{\pi^2}
\int_{q_0(t)}^{q_{1}(t)}
\Re
e\left[(\upsilon^2(t,q)+q^2)\tfpsi(\upsilon(t,q))\frac{\partial\upsilon}{\partial    t}(t,q)\right]dq,
\end{array}
\end{array}\right.$$
$\dsp\hbox{if }t_{0} <t\leq t_{h_2}  \hbox{ and }  \left|\Im
    m\left[\gamma(t_{0},0)\right]\right|
  <\frac{1}{V_{\max}}$,
$$
\left\{
\begin{array}{ll}
  \dsp u^-_{PfS,x}(x,z,t)= -\frac{F^+_{Pf}}{\pi^2}\int_{0}^{q_0(t)}
\Re
e\left[\ic\gamma(t,q)\kapsi(\gamma(t,q))\tfpsi(\gamma(t,q))\frac{\partial
    \gamma}{\partial t}(t,q)\right]dq,
\\[18pt]
\dsp u^-_{PfS,z}(x,z,t)=\frac{\mat
F^+_{Pf}}{\pi^2}\int_{0}^{q_0(t)}
\Re
e\left[(\gamma^2(t,q)+q^2)\tfpsi(\gamma(t,q))\frac{\partial
    \gamma
}{\partial t}(t,q)\right]dq,
\end{array}\right.$$
$\dsp\hbox{if } t_{h_2}<t
\hbox{ and } \left|\Im
    m\left[\gamma(t_{0},0)\right]\right|
  <\frac{1}{V_{\max}}$ or
 $\dsp\hbox{if }t_{0} <t
\hbox{ and } \left|\Im
    m\left[\gamma(t_{0},0)\right]\right|
  \geq\frac{1}{V_{\max}}$\\
and  $\gr
 u^-_{PfS}(x,z,t)=0$ else.\\\\
 $t_0$ denotes here the arrival time of the transmitted $PfS$
 volume wave at point $(x,0,z)$,
\begin{eqnarray}
  t_{h_1}&=&h\sqrt{\frac{1}{{V^+_{Pf}}^2}-\frac{1}{V^2_{\max}}}-z\sqrt{\frac{1}{{V^-_{S}}^2}-\frac{1}{V^2_{\max}}}+\frac{|x|}{V_{\max}}
\end{eqnarray}
denotes the arrival time of the transmitted $PfS$
head wave at point $(x,0,z)$,
\begin{eqnarray}
  t_{h_2}&=&\frac{h^2+z^2-hz\left(\frac{c_2}{c_1}+
\frac{c_1}{c_2}
\right)+x^2}
{\frac{h}{c_1}-\frac{z}{c_2}}
\end{eqnarray}
 denotes the time after which there is no longer head
wave at point $(x,0,z)$, where
$$c_1=\sqrt{\frac{1}{{V_{Pf}^+}^2}-\frac{1}{{V^2_{\max}}}} \hbox{ and } c_2=\sqrt{\frac{1}{{V^-_{S}}^2}-\frac{1}{{V^2_{\max}}}}.$$
The function $q_0:[t_0\,;\,+\infty]\mapsto\setR^+$ is
the reciprocal function of $\tilde t_0
:\setR^+\mapsto:[t_0,+\infty]$, where $\widetilde
t_0(q)$ is the arrival time  at point $(x,0,z)$ of the fictitious 
transmitted $PfS$ volume wave, propagating at a velocity $\mat
V_{Pf}^+(q)$ in the top layer and at velocity $\mat
V^-_{S}(q)$ in the bottom layer.\\[5pt]
The function $q_1:[t_1\,;\,t_0]\mapsto\setR^+$ is defined by 
$$q_{1}(t)=\sqrt{\frac{1}{x^2}\left(t+z\sqrt{\frac{1}{{V^-_{S}}^2}-\frac{1}{V^2_{\max}}}
-h\sqrt{\frac{1}{{V_{Pf}^+}^2}-\frac{1}{V^2_{\max}}}
\right)^2-\frac{1}{V_{\max}^2}}.$$
The function $\gamma:\{(t,q)\in\setR^+\times\setR^+\,|\, t>\tilde t_0(q)\}\mapsto\setC$ is implicitly defined as the only root of the function 
$${\cal
  F}(\gamma,q,t)=-z\left(\frac{1}{{\mat V^-_{S}}^2(q)}+\gamma^2\right)^{1/2}+h\left(\frac{1}{{\mat V_{Pf}^+}^2(q)}+\gamma^2\right)^{1/2}+i\gamma x-t$$
whose real part is positive.\\
The function $\upsilon:E_1\cup E_2\mapsto\setC$ is implicitly defined as the only root of the function 
$${\cal
  F}(\upsilon,q,t)=-z\left(\frac{1}{{\mat V^-_{S}}^2(q)}+\upsilon^2\right)^{1/2}+h\left(\frac{1}{{\mat V_{Pf}^+}^2(q)}+\upsilon^2\right)^{1/2}+i\upsilon x-t$$
such that $\Im m\left[\partial_t \upsilon(t,q)
\right]<0$, with
$$E_1=\left\{
(t,q)\in\setR^+\times\setR^+\,|\, t_{h_1}<t<t_{0} \hbox{ and } 0<q<q_0(t)
\right\} $$ and  $$E_2=\left\{
(t,q)\in\setR^+\times\setR^+\,|\, t_{0}<t<t_{h_1} \hbox{ and } q_0(t)<q<q_1(t)\right\}.$$
  \item  $\gr u^-_{PsPf}$ is the solid displacement of the
    transmitted $PsPf$ wave and satisfies:
$$
\left\{
\begin{array}{ll}
 \dsp u^-_{PsPf,x}(x,z,t)=-\frac{\mat
   P_{11}^-F^+_{Ps}}{\pi^2}\int_{0}^{q_{1}(t)}
\Re
 e\left[\ic\upsilon(t,q)\tsf(\upsilon(t,q))\frac{\partial \upsilon}{\partial t}(t,q)\right]dq,
\\[18pt]
\dsp u^-_{PsPf,z}(x,z,t)= 
\frac{\mat P_{11}^-F^+_{Ps}}{\pi^2}\int_{0}^{q_{1}(t)}
\Re
e\left[\kaf(\upsilon(t,q))\tsf(\upsilon(t,q))\frac{\partial\upsilon}{\partial    t}(t,q)\right]dq,
\end{array}\right.$$
$\dsp\hbox{if }t_{h_1} <t\leq t_{0} \hbox{ and } \left|\Im
    m\left[\gamma(t_{0},0)\right]\right|
  <\frac{1}{V_{\max}},$
$$\left\{
\begin{array}{lcl}
  \dsp u^-_{PsPf,x}(x,z,t)&=&\dsp 
\begin{array}[t]{ll}
-&\dsp\frac{\mat
    P_{11}^-F^+_{Ps}}{\pi^2}
\int_{0}^{q_0(t)}
\Re
e\left[\ic\gamma(t,q)\tsf(\gamma(t,q))\frac{\partial
    \gamma}{\partial t}(t,q)\right]dq
\\[18pt]
-&\dsp\frac{\mat
    P_{11}^-F^+_{Ps}}{\pi^2}
\int_{q_0(t)}^{q_{1}(t)}
\Re
 e\left[\ic\upsilon(t,q)\tsf(\upsilon(t,q))\frac{\partial \upsilon}{\partial t}(t,q)\right]dq,
\end{array}
\\[60pt]
\dsp u^-_{PsPf,z}(x,z,t)&=&\dsp 
\begin{array}[t]{ll}
&\dsp \frac{\mat
  P_{11}^-F^+_{Ps}}{\pi^2}\int_{0}^{q_0(t)}
\Re
e\left[\kaf(\gamma(t,q))\tsf(\gamma(t,q))\frac{\partial
    \gamma
}{\partial t}(t,q)\right]dq
\\[18pt]
+&\dsp\frac{\mat
  P_{11}^-F^+_{Ps}}{\pi^2}
\int_{q_0(t)}^{q_{1}(t)}
\Re
e\left[\kaf(\upsilon(t,q))\tsf(\upsilon(t,q))\frac{\partial\upsilon}{\partial    t}(t,q)\right]dq,
\end{array}
\end{array}\right.$$
$\dsp\hbox{if }t_{0} <t\leq t_{h_2}  \hbox{ and }  \left|\Im
    m\left[\gamma(t_{0},0)\right]\right|
  <\frac{1}{V_{\max}}$,
$$
\left\{
\begin{array}{ll}
  \dsp u^-_{PsPf,x}(x,z,t)=- \frac{\mat P_{11}^-F^+_{Ps}}{\pi^2}\int_{0}^{q_0(t)}
\Re
e\left[\ic\gamma(t,q)\tsf(\gamma(t,q))\frac{\partial
    \gamma}{\partial t}(t,q)\right]dq,
\\[18pt]
\dsp u^-_{PsPf,z}(x,z,t)=\frac{\mat
  P_{11}^-F^+_{Ps}}{\pi^2}\int_{0}^{q_0(t)}
\Re
e\left[\kaf(\gamma(t,q))\tsf(\gamma(t,q))\frac{\partial
    \gamma
}{\partial t}(t,q)\right]dq,
\end{array}\right.$$
$\dsp\hbox{if } t_{h_2}<t
\hbox{ and } \left|\Im
    m\left[\gamma(t_{0},0)\right]\right|
  <\frac{1}{V_{\max}}$ or
 $\dsp\hbox{if }t_{0} <t
\hbox{ and } \left|\Im
    m\left[\gamma(t_{0},0)\right]\right|
  \geq\frac{1}{V_{\max}}$\\
and  $\gr
 u^-_{PsPf}(x,z,t)=0$ else.\\\\
 $t_0$ denotes here the arrival time of the transmitted $PsPf$
 volume wave at point $(x,0,z)$,
\begin{eqnarray}
  t_{h_1}&=&h\sqrt{\frac{1}{{V^+_{Ps}}^2}-\frac{1}{V^2_{\max}}}-z\sqrt{\frac{1}{{V^-_{Pf}}^2}-\frac{1}{V^2_{\max}}}+\frac{|x|}{V_{\max}}
\end{eqnarray}
denotes the arrival time of the transmitted $PsPf$
head wave at point $(x,0,z)$,
\begin{eqnarray}
  t_{h_2}&=&\frac{h^2+z^2-hz\left(\frac{c_2}{c_1}+
\frac{c_1}{c_2}
\right)+x^2}
{\frac{h}{c_1}-\frac{z}{c_2}}
\end{eqnarray}
 denotes the time after which there is no longer head
wave at point $(x,0,z)$, where
$$c_1=\sqrt{\frac{1}{{V_{Ps}^+}^2}-\frac{1}{{V^2_{\max}}}} \hbox{ and } c_2=\sqrt{\frac{1}{{V^-_{Pf}}^2}-\frac{1}{{V^2_{\max}}}}.$$
The function $q_0:[t_0\,;\,+\infty]\mapsto\setR^+$ is
the reciprocal function of $\tilde t_0
:\setR^+\mapsto:[t_0,+\infty]$, where $\widetilde
t_0(q)$ is the arrival time  at point $(x,0,z)$ of the fictitious 
transmitted $PsPf$ volume wave, propagating at a velocity $\mat
V_{Ps}^+(q)$ in the top layer and at velocity $\mat
V^-_{Pf}(q)$ in the bottom layer.\\[5pt]
The function $q_1:[t_1\,;\,t_0]\mapsto\setR^+$ is defined by 
$$q_{1}(t)=\sqrt{\frac{1}{x^2}\left(t+z\sqrt{\frac{1}{{V^-_{Pf}}^2}-\frac{1}{V^2_{\max}}}
-h\sqrt{\frac{1}{{V_{Ps}^+}^2}-\frac{1}{V^2_{\max}}}
\right)^2-\frac{1}{V_{\max}^2}}.$$
The function $\gamma:\{(t,q)\in\setR^+\times\setR^+\,|\, t>\tilde t_0(q)\}\mapsto\setC$ is implicitly defined as the only root of the function 
$${\cal
  F}(\gamma,q,t)=-z\left(\frac{1}{{\mat V^-_{Pf}}^2(q)}+\gamma^2\right)^{1/2}+h\left(\frac{1}{{\mat V_{Ps}^+}^2(q)}+\gamma^2\right)^{1/2}+i\gamma x-t$$
whose real part is positive.\\
The function $\upsilon:E_1\cup E_2\mapsto\setC$ is implicitly defined as the only root of the function 
$${\cal
  F}(\upsilon,q,t)=-z\left(\frac{1}{{\mat V^-_{Pf}}^2(q)}+\upsilon^2\right)^{1/2}+h\left(\frac{1}{{\mat V_{Ps}^+}^2(q)}+\upsilon^2\right)^{1/2}+i\upsilon x-t$$
such that $\Im m\left[\partial_t \upsilon(t,q)
\right]<0$, with
$$E_1=\left\{
(t,q)\in\setR^+\times\setR^+\,|\, t_{h_1}<t<t_{0} \hbox{ and } 0<q<q_0(t)
\right\} $$ and  $$E_2=\left\{
(t,q)\in\setR^+\times\setR^+\,|\, t_{0}<t<t_{h_1} \hbox{ and } q_0(t)<q<q_1(t)\right\}.$$
  \item  $\gr u^-_{PsPs}$ is the solid displacement of the
    transmitted $PsPs$ wave and satisfies:
$$
\left\{
\begin{array}{ll}
 \dsp u^-_{PsPs,x}(x,z,t)=-\frac{\mat
   P_{12}^-F^+_{Ps}}{\pi^2}\int_{0}^{q_{1}(t)}
\Re
 e\left[\ic\upsilon(t,q)\tss(\upsilon(t,q))\frac{\partial \upsilon}{\partial t}(t,q)\right]dq,
\\[18pt]
\dsp u^-_{PsPs,z}(x,z,t)= 
\frac{\mat P_{12}^-F^+_{Ps}}{\pi^2}\int_{0}^{q_{1}(t)}
\Re
e\left[\kas(\upsilon(t,q))\tss(\upsilon(t,q))\frac{\partial\upsilon}{\partial    t}(t,q)\right]dq,
\end{array}\right.$$
$\dsp\hbox{if }t_{h_1} <t\leq t_{0} \hbox{ and } \left|\Im
    m\left[\gamma(t_{0},0)\right]\right|
  <\frac{1}{V_{\max}},$
$$\left\{
\begin{array}{lcl}
  \dsp u^-_{PsPs,x}(x,z,t)&=&\dsp 
\begin{array}[t]{ll}
-&\dsp\frac{\mat
    P_{12}^-F^+_{Ps}}{\pi^2}
\int_{0}^{q_0(t)}
\Re
e\left[\ic\gamma(t,q)\tss(\gamma(t,q))\frac{\partial
    \gamma}{\partial t}(t,q)\right]dq
\\[18pt]
-&\dsp\frac{\mat
    P_{12}^-F^+_{Ps}}{\pi^2}
\int_{q_0(t)}^{q_{1}(t)}
\Re
 e\left[\ic\upsilon(t,q)\tss(\upsilon(t,q))\frac{\partial \upsilon}{\partial t}(t,q)\right]dq,
\end{array}
\\[60pt]
\dsp u^-_{PsPs,z}(x,z,t)&=&\dsp 
\begin{array}[t]{ll}
&\dsp \frac{\mat
  P_{12}^-F^+_{Ps}}{\pi^2}\int_{0}^{q_0(t)}
\Re
e\left[\kas(\gamma(t,q))\tss(\gamma(t,q))\frac{\partial
    \gamma
}{\partial t}(t,q)\right]dq
\\[18pt]
+&\dsp\frac{\mat
  P_{12}^-F^+_{Ps}}{\pi^2}
\int_{q_0(t)}^{q_{1}(t)}
\Re
e\left[\kas(\upsilon(t,q))\tss(\upsilon(t,q))\frac{\partial\upsilon}{\partial    t}(t,q)\right]dq,
\end{array}
\end{array}\right.$$
$\dsp\hbox{if }t_{0} <t\leq t_{h_2}  \hbox{ and }  \left|\Im
    m\left[\gamma(t_{0},0)\right]\right|
  <\frac{1}{V_{\max}}$,
$$
\left\{
\begin{array}{ll}
  \dsp u^-_{PsPs,x}(x,z,t)=- \frac{\mat P_{12}^-F^+_{Ps}}{\pi^2}\int_{0}^{q_0(t)}
\Re
e\left[\ic\gamma(t,q)\tss(\gamma(t,q))\frac{\partial
    \gamma}{\partial t}(t,q)\right]dq,
\\[18pt]
\dsp u^-_{PsPs,z}(x,z,t)=\frac{\mat
  P_{12}^-F^+_{Ps}}{\pi^2}\int_{0}^{q_0(t)}
\Re
e\left[\kas(\gamma(t,q))\tss(\gamma(t,q))\frac{\partial
    \gamma
}{\partial t}(t,q)\right]dq,
\end{array}\right.$$
$\dsp\hbox{if } t_{h_2}<t
\hbox{ and } \left|\Im
    m\left[\gamma(t_{0},0)\right]\right|
  <\frac{1}{V_{\max}}$ or
 $\dsp\hbox{if }t_{0} <t
\hbox{ and } \left|\Im
    m\left[\gamma(t_{0},0)\right]\right|
  \geq\frac{1}{V_{\max}}$\\
and  $\gr
 u^-_{PsPs}(x,z,t)=0$ else.\\\\
 $t_0$ denotes here the arrival time of the transmitted $PsPs$
 volume wave at point $(x,0,z)$,
\begin{eqnarray}
  t_{h_1}&=&h\sqrt{\frac{1}{{V^+_{Ps}}^2}-\frac{1}{V^2_{\max}}}-z\sqrt{\frac{1}{{V^-_{Ps}}^2}-\frac{1}{V^2_{\max}}}+\frac{|x|}{V_{\max}}
\end{eqnarray}
denotes the arrival time of the transmitted $PsPs$
head wave at point $(x,0,z)$,
\begin{eqnarray}
  t_{h_2}&=&\frac{h^2+z^2-hz\left(\frac{c_2}{c_1}+
\frac{c_1}{c_2}
\right)+x^2}
{\frac{h}{c_1}-\frac{z}{c_2}}
\end{eqnarray}
 denotes the time after which there is no longer head
wave at point $(x,0,z)$, where
$$c_1=\sqrt{\frac{1}{{V_{Ps}^+}^2}-\frac{1}{{V^2_{\max}}}} \hbox{ and } c_2=\sqrt{\frac{1}{{V^-_{Ps}}^2}-\frac{1}{{V^2_{\max}}}}.$$
The function $q_0:[t_0\,;\,+\infty]\mapsto\setR^+$ is
the reciprocal function of $\tilde t_0
:\setR^+\mapsto:[t_0,+\infty]$, where $\widetilde
t_0(q)$ is the arrival time  at point $(x,0,z)$ of the fictitious 
transmitted $PsPs$ volume wave, propagating at a velocity $\mat
V_{Ps}^+(q)$ in the top layer and at velocity $\mat
V^-_{Ps}(q)$ in the bottom layer.\\[5pt]
The function $q_1:[t_1\,;\,t_0]\mapsto\setR^+$ is defined by 
$$q_{1}(t)=\sqrt{\frac{1}{x^2}\left(t+z\sqrt{\frac{1}{{\mat V^-_{Ps}}^2}-\frac{1}{V^2_{\max}}}
-h\sqrt{\frac{1}{{\mat V_{Ps}^+}^2}-\frac{1}{V^2_{\max}}}
\right)^2-\frac{1}{V_{\max}^2}}.$$
The function $\gamma:\{(t,q)\in\setR^+\times\setR^+\,|\, t>\tilde t_0(q)\}\mapsto\setC$ is implicitly defined as the only root of the function 
$${\cal
  F}(\gamma,q,t)=-z\left(\frac{1}{{\mat V^-_{Ps}}^2(q)}+\gamma^2\right)^{1/2}+h\left(\frac{1}{{\mat V_{Ps}^+}^2(q)}+\gamma^2\right)^{1/2}+i\gamma x-t$$
whose real part is positive.\\
The function $\upsilon:E_1\cup E_2\mapsto\setC$ is implicitly defined as the only root of the function 
$${\cal
  F}(\upsilon,q,t)=-z\left(\frac{1}{{\mat V^-_{Ps}}^2(q)}+\upsilon^2\right)^{1/2}+h\left(\frac{1}{{\mat V_{Ps}^+}^2(q)}+\upsilon^2\right)^{1/2}+i\upsilon x-t$$
such that $\Im m\left[\partial_t \upsilon(t,q)
\right]<0$, with
$$E_1=\left\{
(t,q)\in\setR^+\times\setR^+\,|\, t_{h_1}<t<t_{0} \hbox{ and } 0<q<q_0(t)
\right\} $$ and  $$E_2=\left\{
(t,q)\in\setR^+\times\setR^+\,|\, t_{0}<t<t_{h_1} \hbox{ and } q_0(t)<q<q_1(t)\right\}.$$
  \item  $\gr u^-_{PsS}$ is the solid displacement of the
    transmitted $PsS$ wave and satisfies:
$$
\left\{
\begin{array}{ll}
 \dsp u^-_{PsS,x}(x,z,t)=-\frac{
   F^+_{Ps}}{\pi^2}\int_{0}^{q_{1}(t)}
\Re
 e\left[\ic\upsilon(t,q)\kapsi(\upsilon(t,q))\tspsi(\upsilon(t,q))\frac{\partial \upsilon}{\partial t}(t,q)\right]dq,
\\[18pt]
\dsp u^-_{PsS,z}(x,z,t)= 
\frac{F^+_{Ps}}{\pi^2}\int_{0}^{q_{1}(t)}
\Re
e\left[(\upsilon^2(t,q)+q^2)\tspsi(\upsilon(t,q))\frac{\partial\upsilon}{\partial    t}(t,q)\right]dq,
\end{array}\right.$$
$\dsp\hbox{if }t_{h_1} <t\leq t_{0} \hbox{ and } \left|\Im
    m\left[\gamma(t_{0},0)\right]\right|
  <\frac{1}{V_{\max}},$
$$\left\{
\begin{array}{lcl}
  \dsp u^-_{PsS,x}(x,z,t)&=&\dsp 
\begin{array}[t]{ll}
-&\dsp\frac{
  F^+_{Ps}}{\pi^2}
\int_{0}^{q_0(t)}
\Re
e\left[\ic\gamma(t,q)\kapsi(\gamma(t,q))\tspsi(\gamma(t,q))\frac{\partial
    \gamma}{\partial t}(t,q)\right]dq
\\[18pt]
-&\dsp\frac{
F^+_{Ps}}{\pi^2}
\int_{q_0(t)}^{q_{1}(t)}
\Re
 e\left[\ic\upsilon(t,q)\kapsi(\upsilon(t,q))\tspsi(\upsilon(t,q))\frac{\partial \upsilon}{\partial t}(t,q)\right]dq,
\end{array}
\\[60pt]
\dsp u^-_{PsS,z}(x,z,t)&=&\dsp 
\begin{array}[t]{ll}
&\dsp \frac{\mat
F^+_{Ps}}{\pi^2}\int_{0}^{q_0(t)}
\Re
e\left[(\gamma^2(t,q)+q^2)\tspsi(\gamma(t,q))\frac{\partial
    \gamma
}{\partial t}(t,q)\right]dq
\\[18pt]
+&\dsp\frac{\mat
 F^+_{Ps}}{\pi^2}
\int_{q_0(t)}^{q_{1}(t)}
\Re
e\left[(\upsilon^2(t,q)+q^2)\tspsi(\upsilon(t,q))\frac{\partial\upsilon}{\partial    t}(t,q)\right]dq,
\end{array}
\end{array}\right.$$
$\dsp\hbox{if }t_{0} <t\leq t_{h_2}  \hbox{ and }  \left|\Im
    m\left[\gamma(t_{0},0)\right]\right|
  <\frac{1}{V_{\max}}$,
$$
\left\{
\begin{array}{ll}
  \dsp u^-_{PsS,x}(x,z,t)= -\frac{F^+_{Ps}}{\pi^2}\int_{0}^{q_0(t)}
\Re
e\left[\ic\gamma(t,q)\kapsi(\gamma(t,q))\tspsi(\gamma(t,q))\frac{\partial
    \gamma}{\partial t}(t,q)\right]dq,
\\[18pt]
\dsp u^-_{PsS,z}(x,z,t)=\frac{\mat
F^+_{Ps}}{\pi^2}\int_{0}^{q_0(t)}
\Re
e\left[(\gamma^2(t,q)+q^2)\tspsi(\gamma(t,q))\frac{\partial
    \gamma
}{\partial t}(t,q)\right]dq,
\end{array}\right.$$
$\dsp\hbox{if } t_{h_2}<t
\hbox{ and } \left|\Im
    m\left[\gamma(t_{0},0)\right]\right|
  <\frac{1}{V_{\max}}$ or
 $\dsp\hbox{if }t_{0} <t
\hbox{ and } \left|\Im
    m\left[\gamma(t_{0},0)\right]\right|
  \geq\frac{1}{V_{\max}}$\\
and  $\gr
 u^-_{PsS}(x,z,t)=0$ else.\\\\
 $t_0$ denotes here the arrival time of the transmitted $PsS$
 volume wave at point $(x,0,z)$,
\begin{eqnarray}
  t_{h_1}&=&h\sqrt{\frac{1}{{V^+_{Ps}}^2}-\frac{1}{V^2_{\max}}}-z\sqrt{\frac{1}{{V^-_{S}}^2}-\frac{1}{V^2_{\max}}}+\frac{|x|}{V_{\max}}
\end{eqnarray}
denotes the arrival time of the transmitted $PsS$
head wave at point $(x,0,z)$,
\begin{eqnarray}
  t_{h_2}&=&\frac{h^2+z^2-hz\left(\frac{c_2}{c_1}+
\frac{c_1}{c_2}
\right)+x^2}
{\frac{h}{c_1}-\frac{z}{c_2}}
\end{eqnarray}
 denotes the time after which there is no longer head
wave at point $(x,0,z)$, where
$$c_1=\sqrt{\frac{1}{{V_{Ps}^+}^2}-\frac{1}{{V^2_{\max}}}} \hbox{ and } c_2=\sqrt{\frac{1}{{V^-_{S}}^2}-\frac{1}{{V^2_{\max}}}}.$$
The function $q_0:[t_0\,;\,+\infty]\mapsto\setR^+$ is
the reciprocal function of $\tilde t_0
:\setR^+\mapsto:[t_0,+\infty]$, where $\widetilde
t_0(q)$ is the arrival time  at point $(x,0,z)$ of the fictitious 
transmitted $PsS$ volume wave, propagating at a velocity $\mat
V_{Ps}^+(q)$ in the top layer and at velocity $\mat
V^-_{S}(q)$ in the bottom layer.\\[5pt]
The function $q_1:[t_1\,;\,t_0]\mapsto\setR^+$ is defined by 
$$q_{1}(t)=\sqrt{\frac{1}{x^2}\left(t+z\sqrt{\frac{1}{{V^-_{S}}^2}-\frac{1}{V^2_{\max}}}
-h\sqrt{\frac{1}{{V_{Ps}^+}^2}-\frac{1}{V^2_{\max}}}
\right)^2-\frac{1}{V_{\max}^2}}.$$
The function $\gamma:\{(t,q)\in\setR^+\times\setR^+\,|\, t>\tilde t_0(q)\}\mapsto\setC$ is implicitly defined as the only root of the function 
$${\cal
  F}(\gamma,q,t)=-z\left(\frac{1}{{\mat V^-_{S}}^2(q)}+\gamma^2\right)^{1/2}+h\left(\frac{1}{{\mat V_{Ps}^+}^2(q)}+\gamma^2\right)^{1/2}+i\gamma x-t$$
whose real part is positive.\\
The function $\upsilon:E_1\cup E_2\mapsto\setC$ is implicitly defined as the only root of the function 
$${\cal
  F}(\upsilon,q,t)=-z\left(\frac{1}{{\mat V^-_{S}}^2(q)}+\upsilon^2\right)^{1/2}+h\left(\frac{1}{{\mat V_{Ps}^+}^2(q)}+\upsilon^2\right)^{1/2}+i\upsilon x-t$$
such that $\Im m\left[\partial_t \upsilon(t,q)
\right]<0$, with
$$E_1=\left\{
(t,q)\in\setR^+\times\setR^+\,|\, t_{h_1}<t<t_{0} \hbox{ and } 0<q<q_0(t)
\right\} $$ and  $$E_2=\left\{
(t,q)\in\setR^+\times\setR^+\,|\, t_{0}<t<t_{h_1} \hbox{ and } q_0(t)<q<q_1(t)\right\}.$$
\end{itemize}
\end{theo}
\begin{remark}
  For the practical computations of the velocities,  we won't have to
  explicitly compute the derivatives of the displacement
  $u$, which would be rather tedious, since
$$\frac{du}{dt}\ast f=u\ast f'. $$
Therefore, we'll only have to compute the derivative of the source function $f$.
\end{remark}

\section{Numerical illustration}
\label{Numerical}
\begin{figure}[h]
      \begin{minipage}{.48\linewidth}
      \centerline{
\setlength{\unitlength}{.9cm}
\begin{picture}(8,6)
        \includegraphics[height=5cm]{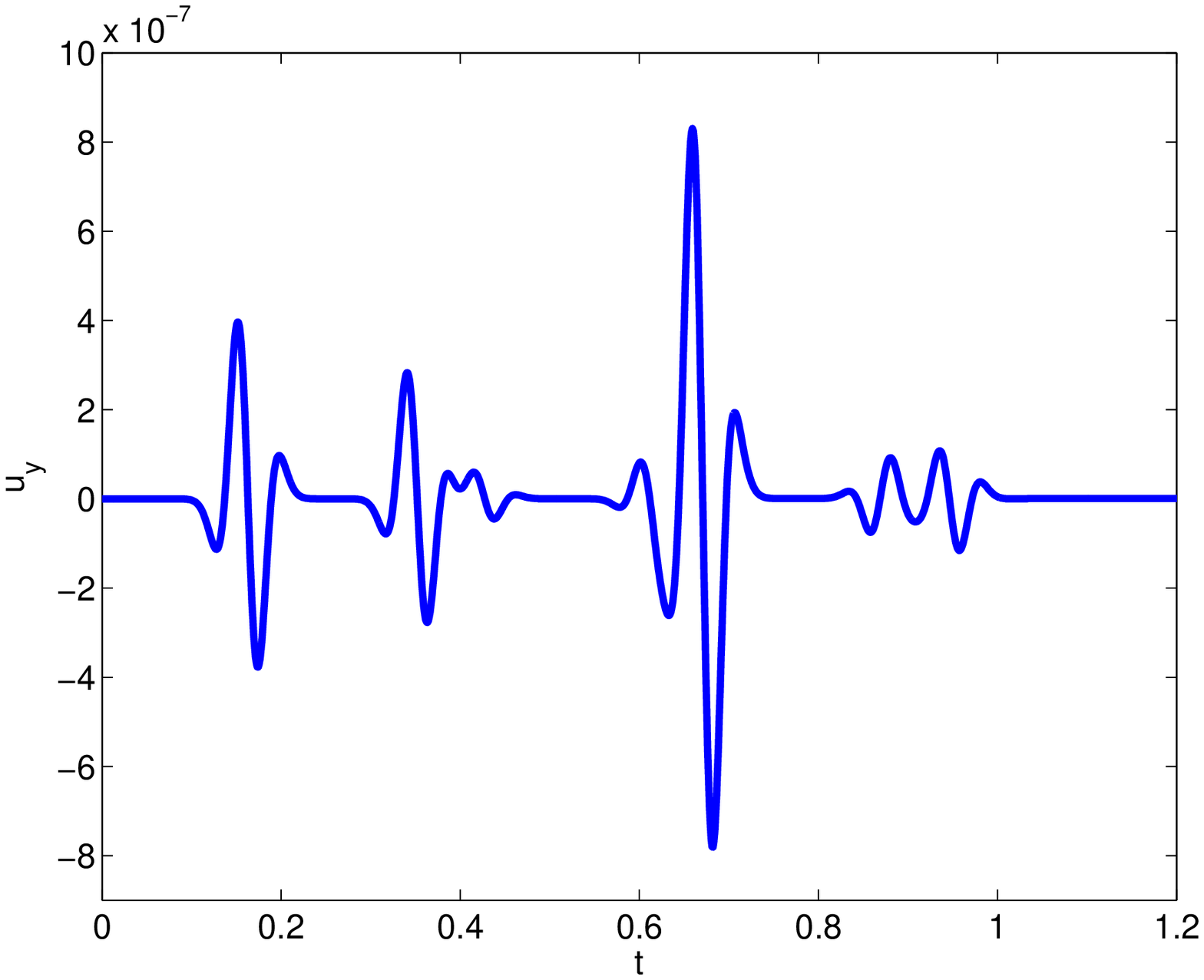}
\psline{<-}(-5,3.4)(-5.2,4)
\put(-5.95,4.5){$Pf$}
\psline{<-}(-4.3,2.25)(-5.1,1.2) 
\put(-6,1.){$Ps$}
\psline{<-}(-3.8,2.4)(-3.9,1.4)
\put(-4.9,1.2){$PfPf$}
\psline{<-}(-3.1,2.5)(-3.5,3)
\put(-4.4,3.4){$PfS$}
\psline{<-}(-2.8,4.2)(-3.3,4.2)
\put(-4.75,4.85){$PsPf$}
\put(-4.35,4.45){\&}
\put(-4.75,4.05){$PfPs$}
\psline{<-}(-1.85,2.7)(-2.1,3.1)
\put(-2.8,3.5){$PsS$}
\psline{<-}(-1.4,2.1)(-1.2,1.3)
\put(-1.8,1.){$PsPs$}
  \end{picture}
      }  
    \end{minipage}
      \begin{minipage}{.48\linewidth}
      \centerline{
\setlength{\unitlength}{.9cm}
\begin{picture}(8,6)
        \includegraphics[height=5cm]{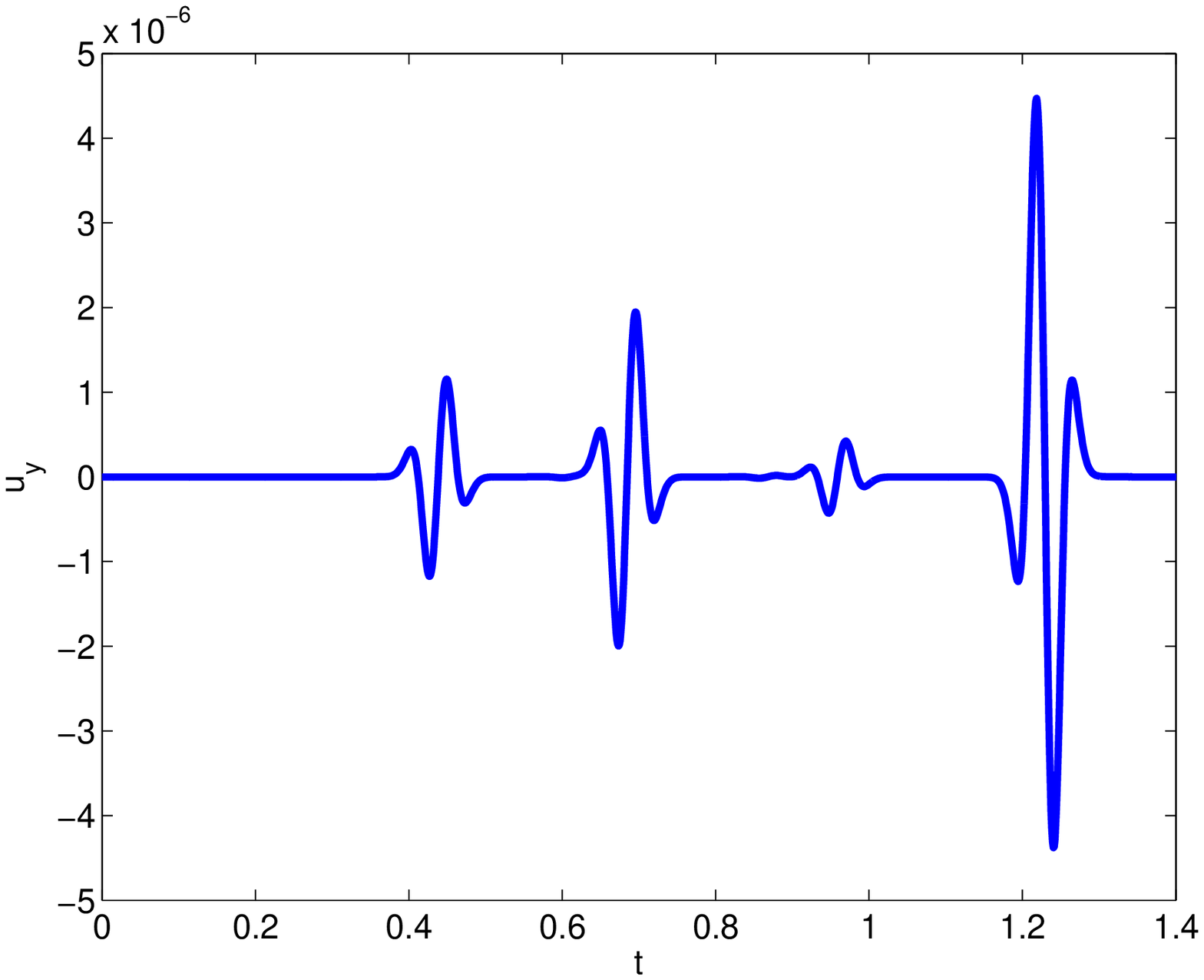}
\psline{<-}(-4.15,2.75)(-5,3.5)
\put(-6.1,4.){$PfPf$}
\psline{<-}(-2.4,2.6)(-2.7,3.6)
\put(-3.5,4.1){$PsS$}
\psline{<-}(-3.4,2.5)(-4,1.5)
\put(-5.,1.3){$PfS$}
\psline{<-}(-3.1,1.6)(-3.1,1.1)
\put(-3.9,.8){$PsPf$}
\psline{<-}(-1.95,2.8)(-1.8,3.6)
\put(-2.5,4.1){$PfPs$}
\psline{<-}(-1.,1.6)(-1.8,1.1)
\put(-2.6,0.8){$PsPs$}
      \end{picture}
    }  
    \end{minipage}
\caption{The $z$ component of the displacement at
  receiver 1 (left picture) and 2 (right
  picture) in the case of a bulk source.}
\label{poroporo3d:fig:1}
  \end{figure}

\begin{figure}[h]
      \begin{minipage}{.48\linewidth}
      \centerline{
\setlength{\unitlength}{.9cm}
\begin{picture}(8,6)
        \includegraphics[height=5cm]{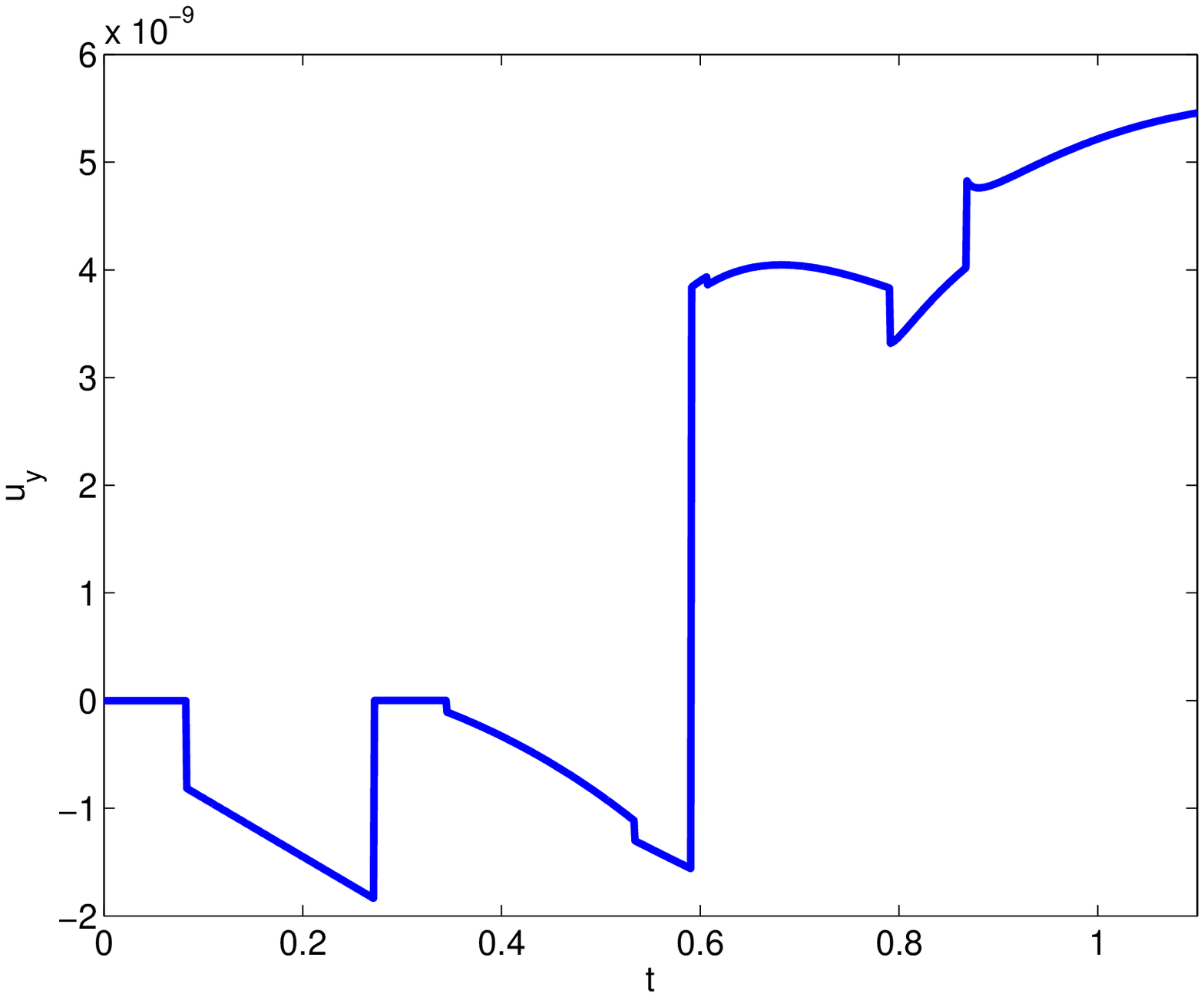}
\psline{<-}(-5.15,1.1)(-4.9,3)
\put(-5.7,3.45){$Pf$}
 \psline{<-}(-4.25,1.55)(-4.5,2.3)
 \put(-5.3,2.6){$Ps$}
 \psline{<-}(-3.9,1.55)(-3.7,2.4)
 \put(-4.7,2.8){$PfPf$}
 \psline{<-}(-2.9,0.85)(-2.3,1.5)
 \put(-2.9,1.8){$PfS$}
 \psline{<-}(-2.7,3.55)(-3.4,4)
 \put(-4.4,4.5){$PsPf$}
 \psline{<-}(-2.55,3.52)(-2.2,3.)
 \put(-2.9,2.9){$PfPs$}
 \psline{<-}(-1.7,3.23)(-1.3,3.0)
 \put(-1.5,3.){$PfS$}
 \psline{<-}(-1.35,4.1)(-1.7,4.3)
 \put(-3,4.8){$PsPs$}
\end{picture}
}  
    \end{minipage}
      \begin{minipage}{.48\linewidth}
      \centerline{
\setlength{\unitlength}{.9cm}
\begin{picture}(8,6)
        \includegraphics[height=5cm]{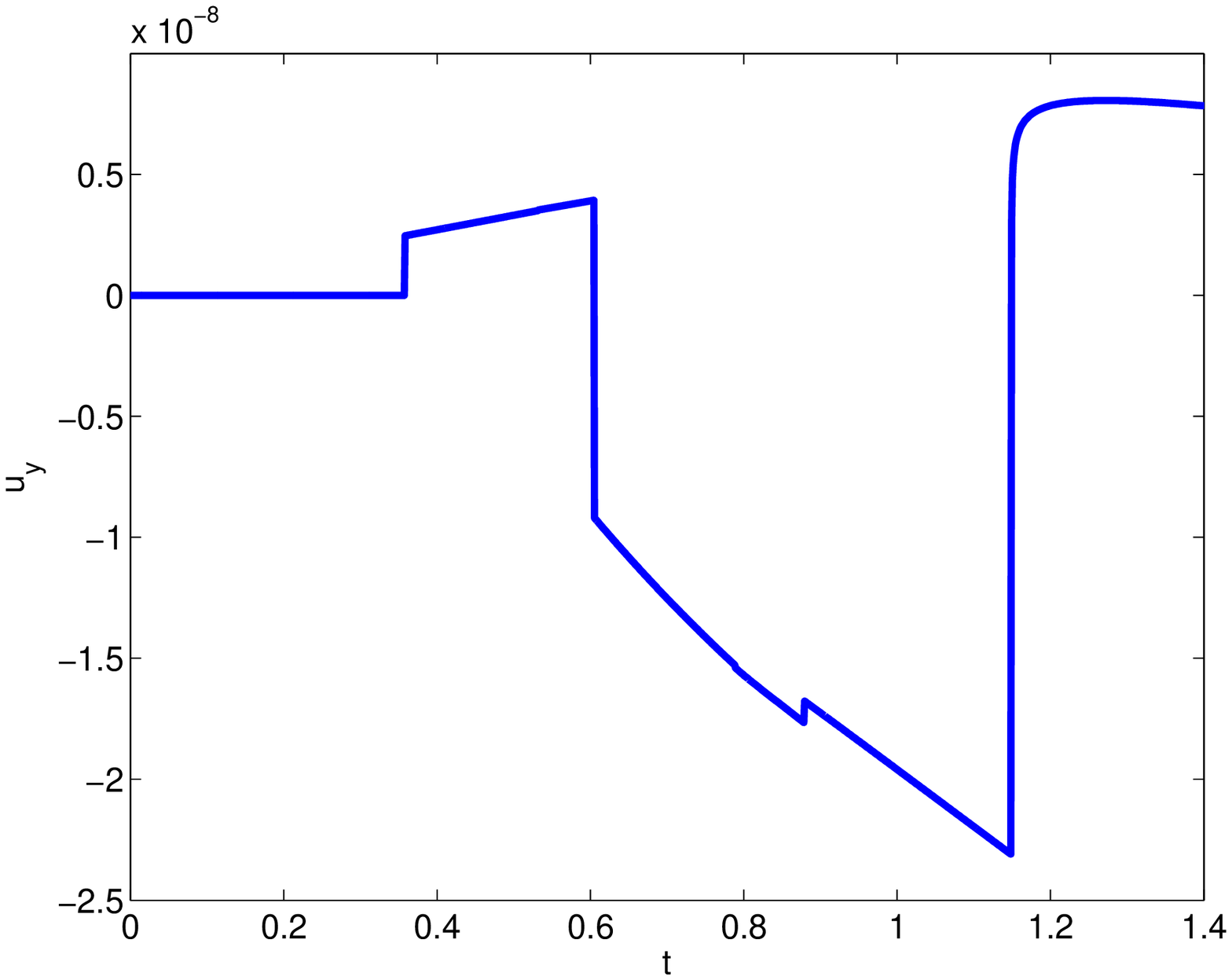}
\psline{<-}(-4.3,3.82)(-5,4.2)
\put(-6.2,4.8){$PfPf$}
\psline{<-}(-3.6,3.95)(-3,4.4)
\put(-3.3,4.8){$PfS$}
 \psline{<-}(-3.4,2.3)(-4.2,1.5)
 \put(-5.65,1.3){$PsPf$}
 \psline{<-}(-2.65,1.7)(-2.5,3.3)
 \put(-3.2,3.8){$PsS$}
 \psline{<-}(-2.3,1.48)(-2.,2)
 \put(-2.8,2.3){$PfPs$}
 \psline{<-}(-1.15,4.35)(-0.8,3)
 \put(-1.4,2.9){$PsPs$}
      \end{picture}
    }  
    \end{minipage}
\caption{The $z$ component of the green function associated to the displacement at
  receiver 1 (left picture) and 2 (right
  picture) in the case of a pressure source.}
\label{poroporo2d:fig:1}
  \end{figure}

\begin{figure}[h]
      \begin{minipage}{.48\linewidth}
      \centerline{
\setlength{\unitlength}{.9cm}
\begin{picture}(8,6)
        \includegraphics[height=5cm]{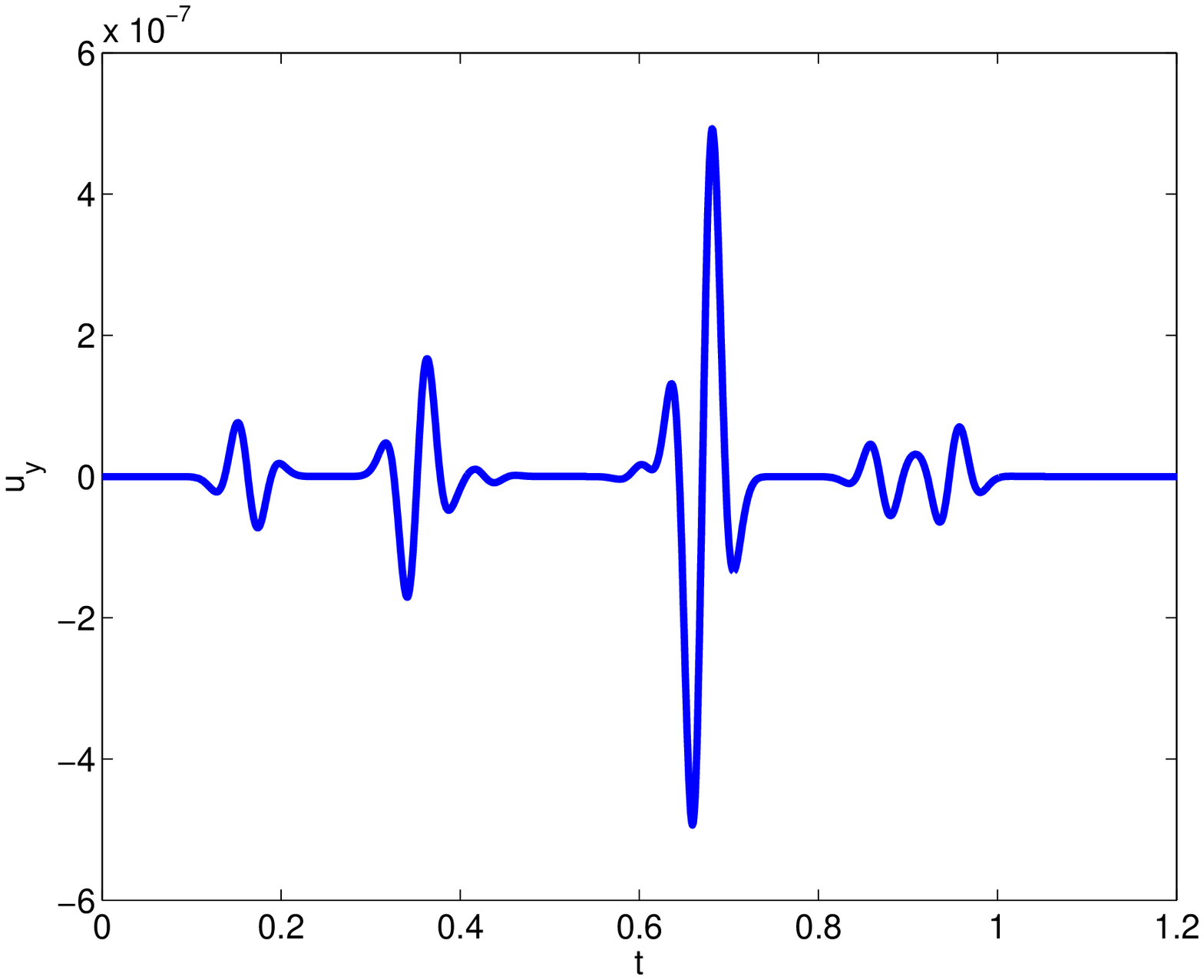}
\psline{<-}(-5,2.9)(-5.2,3.5)
\put(-6.05,4.){$Pf$}
\psline{<-}(-4.2,2.)(-5.1,1.2) 
\put(-6,1.){$Ps$}
\psline{<-}(-3.8,2.5)(-3.9,1.4)
\put(-4.9,1.2){$PfPf$}
\psline{<-}(-3.05,2.65)(-3.5,3)
\put(-4.4,3.4){$PfS$}
\psline{<-}(-2.7,4.2)(-3.3,4.2)
\put(-4.75,4.85){$PsPf$}
\put(-4.35,4.45){\&}
\put(-4.75,4.05){$PfPs$}
\psline{<-}(-1.85,2.8)(-2.1,3.1)
\put(-2.8,3.5){$PsS$}
\psline{<-}(-1.4,2.3)(-1.2,1.5)
\put(-1.8,1.3){$PsPs$}
\end{picture}
}  
    \end{minipage}
      \begin{minipage}{.48\linewidth}
      \centerline{
\setlength{\unitlength}{.9cm}
\begin{picture}(8,6)
        \includegraphics[height=5cm]{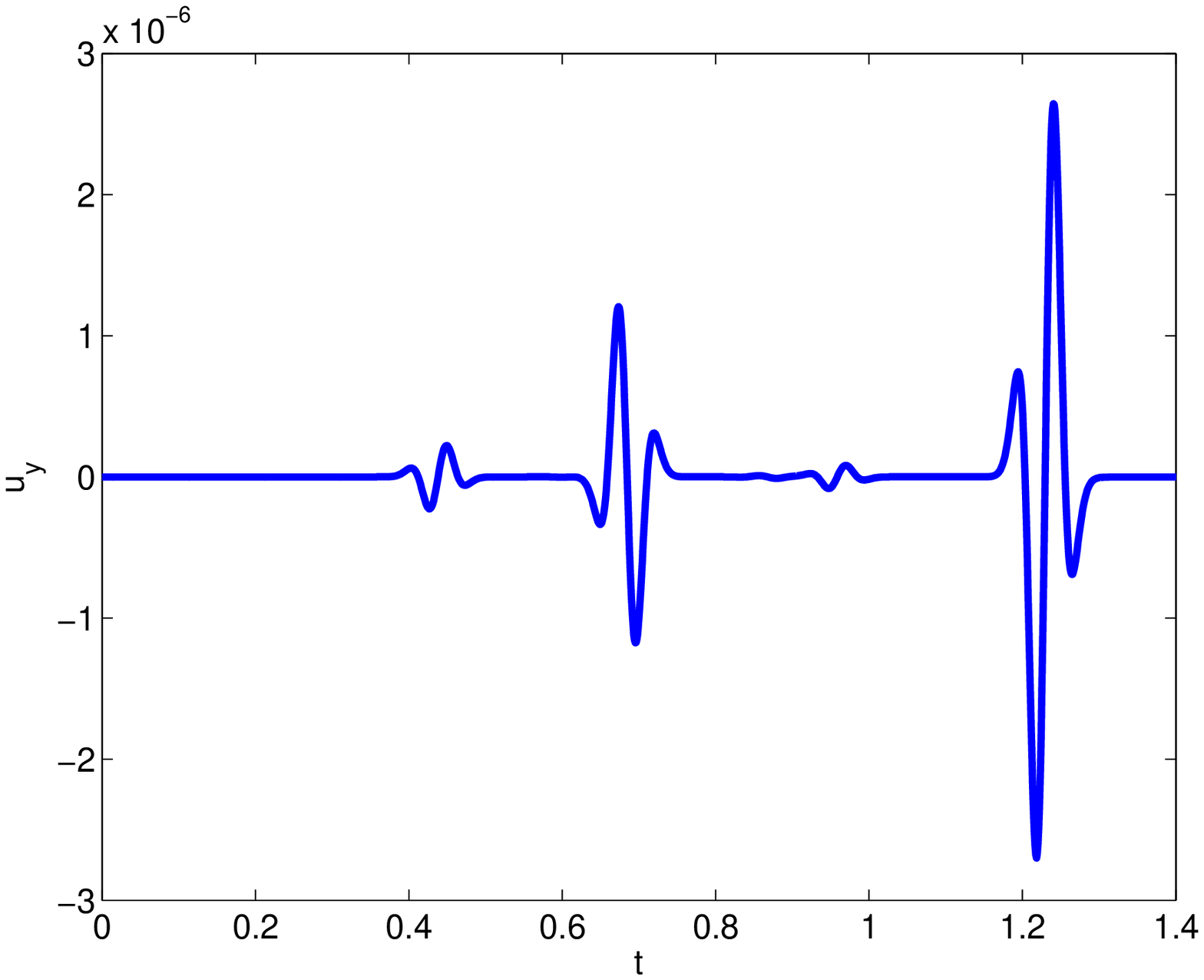}
\psline{<-}(-4.1,2.65)(-5,3.2)
\put(-6.2,3.6){$PfPf$}
\psline{<-}(-3.05,3.45)(-2.8,4)
\put(-3.65,4.55){$PsPf$}
\psline{<-}(-2.3,2.6)(-2.,3.4)
\put(-2.5,3.9){$PsS$}
\psline{<-}(-3.4,2.5)(-3.6,1.5)
\put(-4.6,1.3){$PfS$}
\psline{<-}(-1.95,2.5)(-1.75,2.)
\put(-2.5,1.8){$PfPs$}
\psline{<-}(-1.1,1.)(-1.5,1.)
\put(-2.8,1.){$PsPs$}
      \end{picture}
    }  
    \end{minipage}
\caption{The $z$ component of the displacement at
  receiver 1 (left picture) and 2 (right
  picture) in the case of a pressure source.}
\label{poroporo3d:fig:2}
  \end{figure}

To illustrate our results, we have
computed the green function and the analytical solution to the following problem: we consider an two-layered poroelastic medium whose characteristic coefficients are
\begin{itemize}
\item the solid density: $\rho_s^+=\unit{2200}{\kilo\gram\per\cubic\metre}$ and $\rho_s^-=\unit{2650}{\kilo\gram\per\cubic\metre}$;
\item the fluid density: $\rho_f^+=\unit{950}{\kilo\gram\per\cubic\metre}$  and $\rho_f^-=\unit{750}{\kilo\gram\per\cubic\metre}$ ;
\item the porosity: $\phi^+=0.4$ and $\phi^-=0.2$ ;
\item the tortuosity: $a^+=2$ and $a^-=2$;
\item the solid bulk modulus: $K^+_s=\unit{6.9}{\giga\pascal}$ and $K^-_s=\unit{37}{\giga\pascal}$;
\item the fluid bulk modulus: $K^+_f=\unit{2}{\giga\pascal}$ and $K^-_f=\unit{1.7}{\giga\pascal}$;
\item the frame bulk modulus: $K^+_b=\unit{6.7}{\giga\pascal}$ and  $K^-_b=\unit{2.2}{\giga\pascal}$;
\item the frame shear modulus $\mu^+=\unit{3}{\giga\pascal}$ and $\mu^-=\unit{4.4}{\giga\pascal}$;
\end{itemize}
so that the celerity of the waves in the
poroelastic medium are:
\begin{itemize}
\item for the fast P wave, $\Vpf^+=\unit{2692}{\meter\per\second}$ and  $\Vpf^-=\unit{2535}{\meter\per\second}$;
\item for the slow P wave, $\Vps^+=\unit{1186}{\meter\per\second}$ and  $\Vps^-=\unit{744}{\meter\per\second}$;
\item for the $\psi$ wave, $\Vs^+=\unit{1409}{\meter\per\second}$ and $\Vs^-=\unit{1415}{\meter\per\second}.$
\end{itemize}
The source is located in the top layer, at
\unit{500}{\meter} from the interface. 
We used two types of sources in space: the first one is a bulk source such that 
$f_u=f_w=-10^{10}$ and $f_p=0$; the second one is a pressure source such that $f_u=f_w=0$ and $f_p=1.$ In each case we used a fourth derivative of a
Gaussian of dominant frequency
$f_0=\unit{15}{\hertz}$:
$$f(t)=2\frac{\pi^2}{f_0^2}\left[3+12\frac{\pi^2}{f_0^2}\left(t-\frac{1}{f_0}\right)^2+4
  \frac{\pi^4}{f_0^4}\left(t-\frac{1}{f_0}\right)^4\right]e^{-\frac{\pi^2}{f_0^2}\left(t-\frac{1}{f_0}\right)^2}$$
for the source in time.  We compute the solution
at two receivers, the first one is in the upper
layer, at \unit{533}{\meter} from the interface;
the second one is in the bottom layer, at
\unit{533}{\meter} from the interface; both are
located on a vertical line at \unit{400}{\meter}
from the source (see
Fig.~\ref{validation2d:fig:6}).  We represent the
$z$ component of the green function associated to the solid displacement from $t=0$ to
$t=\unit{1.4}{\second}$ in
Fig.~\ref{acousporo:fig:2} for the bulk source and
in Fig.~\ref{poroporo2d:fig:1} for the pressure source. 
In Figs.~\ref{poroporo3d:fig:1} and~\ref{poroporo3d:fig:2}, we plot  the solid displacement.
The left pictures represents the solution
at receiver 1 while the right pictures represents
the solution at receiver 2. As all the types of waves are computed independently, it is easy to distinguish all of them, as it is indicated in the figures. 
solution. \\[10pt]
\begin{figure}[H]
  \centering
\setlength{\unitlength}{.9cm}
  \begin{picture}(6,6.5)(2,0)
\psset{xunit=1.5cm}
\psframe(0,.25)(6,5.75)
\psline(0,3)(6,3)
\put(1.5,4.3){$\Omega^+$}
\put(1.5,2.1){$\Omega^-$}
\put(4.3,5.2){Source}
\psdot(3,4.5)
\psline{<->}(3,3)(3,4.45)
\put(3.5,4.){\unit{500}{\metre}}
\psline{<->}(5,3)(5,4.95)
\put(8.5,4.3){\unit{533}{\metre}}
\psline{<->}(5,3)(5,1.05)
\put(8.5,2){\unit{533}{\metre}}
\psline{<->}(3.05,4.5)(5,4.5)
\put(5.9,4.6){\unit{400}{\metre}}
\psdot(5,5)
\psdot(5,1)
\put(7.4,5.8){Receiver 1}
\put(7.4,0.7){Receiver 2}
  \end{picture}
  \caption{Configuration of the experiment}
\label{validation2d:fig:6}
\end{figure}
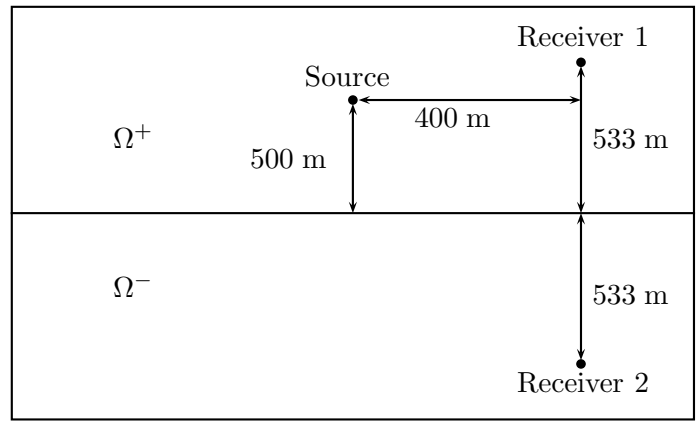
    \begin{figure}[h]
      \begin{minipage}{.48\linewidth}
      \centerline{
\setlength{\unitlength}{.9cm}
\begin{picture}(8,6)
        \includegraphics[height=5cm]{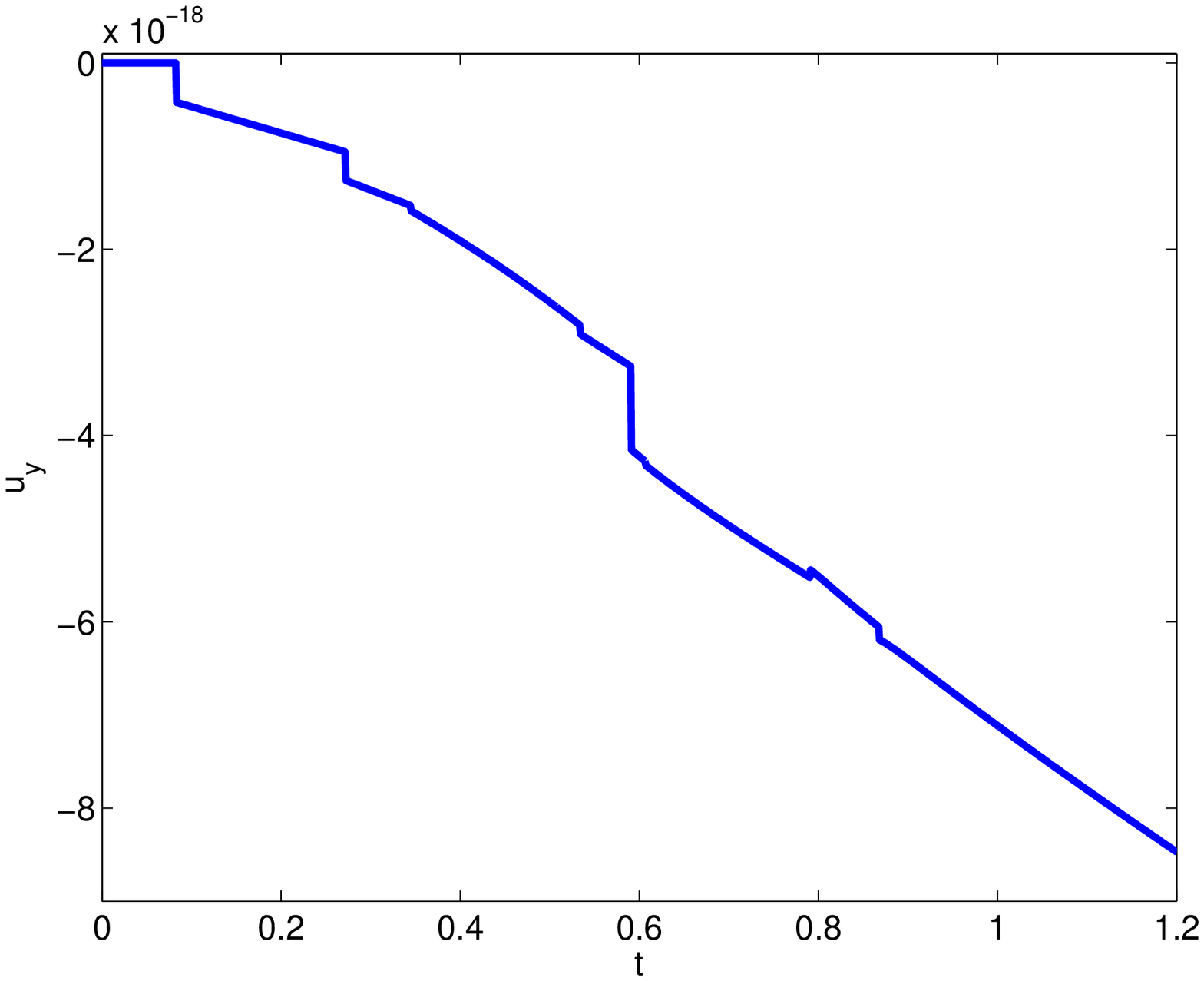}
\psline{<-}(-5.3,4.35)(-5.4,3.8)
\put(-6.3,3.8){$Pf$}
\psline{<-}(-4.45,4.0)(-4.8,3)
\put(-5.7,3.){$Ps$}
\psline{<-}(-4.1,3.85)(-4.6,2)
\put(-5.7,1.8){$PfPf$}
\psline{<-}(-3.3,3.2)(-3.7,2.85)
\put(-4.7,2.85){$PfS$}
\psline{<-}(-3,2.75)(-2.2,4.2)
\put(-2.9,4.7){$PsPf$}
\psline{<-}(-3,2.6)(-3.3,2.2)
\put(-4.4,2.){$PfPs$}
\psline{<-}(-2.17,2.)(-2.5,1.6)
\put(-3.3,1.5){$PsS$}
\psline{<-}(-1.75,1.8)(-1.3,2.4)
\put(-1.9,2.7){$PsPs$}
  \end{picture}
      }  
    \end{minipage}
      \begin{minipage}{.48\linewidth}
      \centerline{ 
\setlength{\unitlength}{.9cm}
\begin{picture}(8,6)
        \includegraphics[height=5cm]{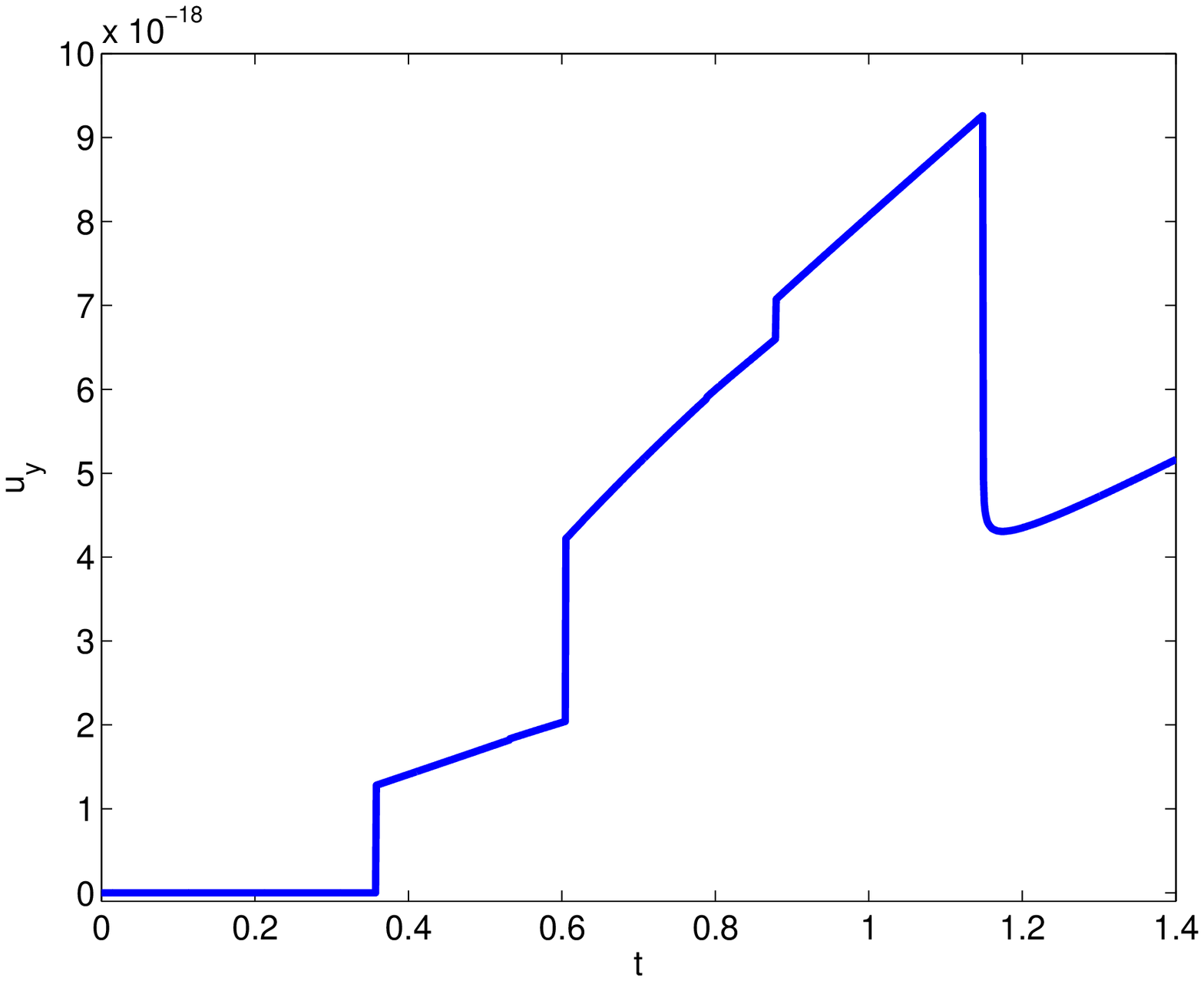}
\psline{<-}(-4.3,1)(-4.7,1.4)
\put(-5.9,1.6){$PfPf$}
\psline{<-}(-3.6,1.3)(-4,1.7)
\put(-5.,2.){$PfS$}
\psline{<-}(-3.4,2.3)(-3.8,2.8)
\put(-4.8,3.2){$PsPf$}
\psline{<-}(-2.3,3.5)(-2.8,4)
\put(-3.6,4.5){$PfPs$}
\psline{<-}(-2.7,3)(-3.2,3.5)
\put(-4,4){$PsS$}
\psline{<-}(-1.3,2.3)(-1.8,1.9)
\put(-2.55,1.8){$PsPs$}
      \end{picture}
    }  
    \end{minipage}
\caption{The $z$ component of the green function associated to the displacement at
  receiver 1 (left picture) and 2 (right
  picture), in the case of a bulk source }
\label{acousporo:fig:2}
  \end{figure}
\bibliography{de-rr}
\bibliographystyle{plain}
\tableofcontents

\end{document}